\documentstyle[11pt]{article}

\begin{document}
\normalbaselineskip=18pt
\baselineskip=\normalbaselineskip

%

\renewcommand{\labelenumi}{(\theenumi)}
\renewcommand{\theenumi}{\roman{enumi}}
\renewcommand{\theenumii}{\arabic{enumii}}
\renewcommand{\theenumiii}{\Alph{enumiii}}
\renewcommand{\theenumiv}{\alph{enumiv}}

%
%
%
\newtheorem{lemma}{Lemma}
\makeatletter
\def\thelemma{\arabic{section}.\arabic{lemma}}
\@addtoreset{lemma}{section}
\makeatother
\newtheorem{th}[lemma]{Theorem}
\newtheorem{rem}[lemma]{Remark}
\newtheorem{df}[lemma]{Definition}
\newtheorem{prop}[lemma]{Proposition}
\newtheorem{cor}[lemma]{Corollary}
\newtheorem{conj}[lemma]{Conjecture}
\newtheorem{corollary}[lemma]{Corollary}
\newtheorem{conjecture}[lemma]{Conjecture}
\newtheorem{definition}[lemma]{Definition}
\newtheorem{example}[lemma]{Example}
\newtheorem{problem}[lemma]{Problem}
%
%
%

\newcommand{\qed}{\hfill\hbox{\rule{3pt}{6pt}}}
\def\ns{\negthinspace } 
\def\implies{$\Longrightarrow\ $}
\def\iff{$\Longleftrightarrow\ $}
\def\dist{\partial}
\def\dim{{\rm dim}\;}
\def\Span{{\rm Span}}
\def\proof{\medskip\noindent {\it Proof.}\quad }
\def\bs{\backslash}
\def\eps{\varepsilon}

%
%
\makeatletter
\font\tenbbb=msbm10 at 12pt
\font\sevenbbb=msbm10 
\font\fivebbb=msbm8
\newfam\bbbfam
\textfont\bbbfam=\tenbbb
\scriptfont\bbbfam=\sevenbbb
\scriptscriptfont\bbbfam=\fivebbb
\edef\bbbfam@{\hexnumber@\bbbfam}
\def\Bbb#1{{\fam\bbbfam\relax#1}}
\makeatother
%

\def\R{{\Bbb R}} \def\rr{{\Bbb R}} 
\def\C{{\Bbb C}} 
\def\Z{{\Bbb Z}} 
\def\openC{{\bf\rm C}\kern-.4em 
	   {\vrule height 1.4ex width.08em depth-.04ex}}
\def\we{{\wedge }}
\def\ve{{\vee }}
\def\rk{{\hbox{rank}}}
\def\e{{\hbox{E}}}
\def\l{{\qquad \over \qquad}}
\def\di{{\cal D}}
\def\ma{{\hbox{max}}}
\def\mi{{\hbox{min}}}
\def\p{{\cal P}}
\def \g{{\Gamma = (X, R)}}
\def\de{\bf}

%


\pagestyle{empty} \vspace{4cm} {\LARGE
\begin{center}
           TIGHT DISTANCE-REGULAR GRAPHS
\end{center}} 

\vskip 5em

\begin{center}%
\begin{minipage}[t]{8cm}%
\begin{center}%
{\large Aleksandar Juri\v{s}i\'{c}\\%
Department of Combinatorics and Optimization\\%
University of Waterloo}\\%
\end{center}%
\end{minipage}
\begin{minipage}[t]{8cm}%
\begin{center}%
{\large Jack Koolen\\
Department of Mathematics\\%
University of Eindhoven}\\%
\end{center}%
\end{minipage}%
\end{center}%
\vskip 2em%
\begin{center}%
{\large Paul Terwilliger\\%
Department of Mathematics\\%
University of Wisconsin}\\%
\vskip 1.5mm%
\end{center}

\vskip 1.5em

\begin{center}
                  {\large June 28, 1998}
\end{center}

\bigskip 
\begin{abstract} \noindent
We consider a distance-regular graph $\Gamma $  with diameter $d \ge 3$
and eigenvalues $k=\theta_0>\theta_1>\cdots >\theta_d$. 
We show the intersection numbers $a_1, b_1$ satisfy 
$$
\left(\theta_1 + {k \over a_1+1}\right)
\left(\theta_d + {k \over a_1+1}\right)
\ge - {ka_1b_1 \over (a_1+1)^2}.                    
$$
We say $\Gamma $ is {\it tight} whenever $\Gamma $ is not bipartite, 
and equality holds above. We characterize the tight property in a 
number of ways. 
For example, we show $\Gamma $ is tight if and only if the intersection 
numbers are given by certain rational expressions  involving $d$ 
independent parameters. 
We  show $\Gamma$ is tight if and only if  $a_1\not=0$, $a_d=0$, and
$\Gamma $ is  1-homogeneous in the sense of Nomura.
We show $\Gamma $ is tight if and only if each local graph is connected
strongly-regular, with nontrivial eigenvalues $-1-b_1(1+\theta_1)^{-1}$
and  $-1-b_1(1+\theta_d)^{-1}$. 
Three infinite families and nine sporadic examples of tight 
distance-regular graphs are given.
\end{abstract}

%
%

\bigskip 
\def\leaderfill{\leaders\hbox to 1em{\hss.\hss}\hfill}
{\narrower {\narrower {\narrower {\narrower \small
                   \centerline{CONTENTS}

\medskip\noindent
1. Introduction                                        \leaderfill 1  \\
2. Preliminaries                                       \leaderfill 2 \\
3. Edges that are tight with respect to an eigenvalue  \leaderfill 7 \\
4. Tight edges and combinatorial  regularity           \leaderfill 9 \\
5. The tightness of an edge                           \leaderfill 11 \\
6. Tight graphs and the Fundamental Bound             \leaderfill 12 \\
7. Two characterizations of tight graphs              \leaderfill 14 \\
8. The auxiliary parameter                            \leaderfill 16 \\
9. Feasibility                                        \leaderfill 17 \\
10. A parametrization                                 \leaderfill 19 \\
11. The 1-homogeneous property                        \leaderfill 23 \\
12. The local graph                                   \leaderfill 28 \\
13. Examples of tight distance-regular graphs         \leaderfill 30 \\
References                                            \leaderfill 33 \\

\baselineskip=\normalbaselineskip \par} \par} \par} \par}

\newpage \pagestyle{plain} \setcounter{page}{1}

\section{Introduction} 

Let $\Gamma=(X,R)$ denote a  distance-regular graph with diameter 
$d \geq 3$,  and  eigenvalues 
$\;k=\theta_0 > \theta_1 > \cdots > \theta_d$ (see Section 2 for
definitions). We show the intersection numbers $a_1, b_1$
satisfy
\begin{eqnarray}
 \left(\theta_1 + {k \over a_1+1}\right)\left(\theta_d + {k \over a_1+1}
\right)
\ge - {ka_1b_1 \over (a_1+1)^2}.                       \label{introA}
\end{eqnarray}
We define $\Gamma $ to be {\it tight} whenever $\Gamma $ is not 
bipartite, and equality holds in (\ref{introA}). 
We characterize the tight condition in the following  ways.

\medskip
\noindent  Our first characterization is linear algebraic. 
For all vertices $x\in X$, let $\hat x$ denote the vector in 
$\R^X$ with a $1$ in coordinate $x$, and $0$ in all other coordinates.
Suppose for the moment that
$a_1\not=0$, let  $x, y$ denote adjacent vertices in  $X$, and write
$w= \sum {\hat z}$, where the sum is over all
vertices $z\in X$ adjacent to  both $x$ and $y$.
Let $\theta $ denote one of $\theta_1, \theta_2, \ldots,\theta_d$,
and let
$E$ denote the corresponding primitive idempotent of the Bose-Mesner
algebra. We say
the edge $xy$ is {\it tight with respect to} $\theta $ whenever
$E \hat x$, 
$E \hat y$, 
$E  w$ are linearly dependent.
We show that if
$xy$ is tight with respect to $\theta$,
then $\theta $ is one of $\theta_1, \theta_d$.
Moreover, we show the following are equivalent:
(i)
$\Gamma $ is tight; (ii) $a_1\not=0$ and  
all  edges of $\Gamma $  are   tight with respect to both
$\theta_1, \theta_d$; (iii) 
 $a_1\not=0$ and there exists an  
edge of $\Gamma $  which is  tight with
respect to both  
$\theta_1, \theta_d$. 

\medskip
\noindent Our second characterization of the tight condition
involves the
intersection numbers.
We show  $\Gamma $ is tight
if and only if the intersection numbers
are given by certain rational  expressions involving $d$ independent
variables. 

\medskip
\noindent Our third characterization of the tight condition involves
the concept of  {\it 1-homogeneous}  that appears in the work of  Nomura
\cite{Nom1}, \cite{Nom2}, \cite{Nom3}. See also  Curtin \cite{Cur}.
We show the following are equivalent:
(i) $\Gamma $ is tight;
 (ii) $a_1\not=0, a_d=0$, and $\Gamma $ is 1-homogeneous;
(iii) $a_1\not=0, a_d=0$, and $\Gamma $ is 1-homogeneous
with respect to at least one edge.

\medskip
\noindent 
Our fourth characterization of the tight condition
involves the local structure and is reminiscent of some results 
by Cameron, Goethals and Seidel \cite{CGS} and Dickie and 
Terwilliger \cite{DT}.
For all  $x \in X$, let  $\Delta(x)$ denote the vertex subgraph of 
$\Gamma$ induced on the vertices in $X$ adjacent to $x$.
For notational convenience, define 
 $b^+:=-1-b_1(1+\theta_d)^{-1}$ and 
 $b^-:=-1-b_1(1+\theta_1)^{-1}$. 
We show  the following
are equivalent: (i) $\Gamma$  is tight; (ii) for all $x\in X$,
 $\Delta(x)$ is connected strongly-regular with 
nontrivial eigenvalues $b^+$, $b^-$;
 (iii)  there exists $ x\in X$ such that $\Delta(x)$ is
 connected strongly-regular with 
nontrivial eigenvalues $b^+$, $b^-$.

\medskip
\noindent 
%
We present three infinite families and nine sporadic examples of tight 
distance-regular graphs. These are the Johnson graphs $J(2d,d)$, 
the halved cubes $ \frac12  H(2d,2)$, the Taylor graphs~\cite{Ta},
%
%
four 3-fold antipodal covers of diameter 4 constructed from 
the sporadic Fisher groups~\cite[p. 397]{BCN},
%
%
two 3-fold antipodal covers of diameter 4 constructed by 
Soicher~\cite{Soi},
a 2-fold and a 4-fold antipodal cover of diameter 4 constructed by
Meixner~\cite{Mei}, and the Patterson graph \cite[Thm. 13.7.1]{BCN}, 
which is primitive, distance-transitive and of diameter 4.

\medskip

\baselineskip=\normalbaselineskip


\section{Preliminaries}            

In this section,  we review some definitions and basic concepts. 
See the books of Bannai and Ito \cite{BI} or Brouwer, Cohen, and 
Neumaier \cite{BCN} for more background information. 

\noindent 
Let $\Gamma=(X,R)$ denote a finite, undirected, connected graph,
without loops or multiple edges, with vertex set $X$, edge set
$R$, path-length distance function  $\partial$, and diameter
$d := \hbox{max}\lbrace \partial(x,y)\;\vert\;x,y\in X\rbrace.$
For all $x \in X$ and for all integers $i$, we set 
$\Gamma_i(x):= \lbrace y\in X\;\vert\;\partial(x,y)=i\rbrace.$
We abbreviate $\Gamma(x):=\Gamma_1(x)$.
By the {\it valency} of a vertex  $x\in X$, we mean the 
cardinality of $\Gamma(x)$. Let $k$ denote a nonnegative integer.
Then $\Gamma $ is said to be {\it regular, with valency k}, whenever
each vertex in  $X$ has valency $k$. $\Gamma$ is said to be 
{\it distance-regular} whenever for all integers
$\;h,\;i,\;j$ $\;(0 \le h,\;i,\;j\leq d)$, and for all $\;x,\;y \in X\;$
with $\;\partial (x,y) = h$, the number
$$
p^h_{ij} := \vert \Gamma_i(x) \cap \Gamma_j(y)\vert  
$$
is independent of $\;x$ and $\;y$. The constants $\;p^h_{ij}\;$
are known as the {\it intersection numbers} of $\;\Gamma$. 
For notational convenience,
set $\;c_i:= p^i_{1\,i-1}$ $\;(1 \le i \leq d)$, $\;a_i := p^i_{1i}$ $\;
(0 \leq i \leq d)$, $\;b_i := p^i_{1\,i+1}$ $\;(0 \le i \le d-1)$, $\;
k_i := p^0_{ii}$ $(0\leq i \leq d)$,  and  define
$\;c_0 = 0$, $\;b_d = 0$. We note $a_0=0$ and   $c_1=1$.

\medskip \noindent 
From now on,  $\;\Gamma=(X,R) \;$ will denote a  distance-regular graph 
of diameter $d\geq 3$. Observe  $\Gamma$ is regular with valency 
$k=k_1=b_0$, and that  
\begin{equation}
\quad  k\;=\; c_i \;+\;a_i\;+\;b_i \qquad \qquad \qquad(0 \le i \le d). 
\label{pre2}
\end{equation}
We now recall the Bose-Mesner algebra.
Let $\hbox{Mat}_X(\R)$ 
denote the $\R$-algebra  consisting of all  matrices with entries in
$\R$ whose rows and columns are indexed by $X$.  
For each integer $\;i$ $\;(0 \le i \le d)$, let $A_i$ denote the matrix
in $\hbox{Mat}_X(\R)$ with
%
$\;x,\,y\;$
entry
$$
\bigl( A_i\bigr)_{xy}\;=\;\cases{1, & if $\;\partial(x,y) = i$,\cr
0, & if $\;\partial(x,y) \ne i$\cr}\qquad \qquad (x,\;y \in X).
$$
$A_i$ is known as the {\it ith distance matrix} of $\Gamma $. 
Observe 
\begin{eqnarray}
 &&A_0 = I, \label{eq:a0}
 \\ 
&&A_0+A_1+\ldots + A_d =  J \qquad (J = \hbox{all 1's matrix} ), 
                                                     \label{eq:asum}\\ 
&&A^t_i = A_i \qquad\qquad \qquad \qquad \qquad (0 \le i \le d),
                                                     \label{eq:aitran}\\
&&A_iA_j = \sum^d_{h=0} p^h_{ij} A_h \qquad \quad\qquad 
(0 \le i,\;j\le d).                                \label{eq:aiprod}
\end{eqnarray}
We abbreviate $A:=A_1$, and refer to this as the
{\it adjacency matrix} of $\Gamma$.
Let $M$ denote the subalgebra of  
$\hbox{Mat}_X(\R)$ generated by $A$. We refer to $M$ as the
{\it Bose-Mesner algebra} of $\Gamma $.
Using (\ref{eq:a0})--(\ref{eq:aiprod}), one can readily show
 $\;A_0,\,A_1,\ldots,\,A_d\;$ form a basis
for $M$.
By \cite[p59, p64]{BI}, the algebra $\;M\;$ has a second basis 
$\;E_0,\;E_1,\,\ldots,\;E_d\;$ such that
\begin{eqnarray}
&&E_0 = \vert X \vert^{-1} J, \label{eq:e0} 
\\
&& E_0 + E_1 + \ldots + E_d = I, \label{eq:eisum}
\\
&&E^t_i = E_i \qquad \qquad \qquad \qquad \qquad (0 \le i \le d),
\label{eq:eitrans}
\\
&&E_iE_j =  \delta_{ij}E_i \qquad \qquad \qquad \qquad (0 \le i,\;j\le d).
\label{eq:eimult}
\end{eqnarray}
The $\;E_0,\;E_1,\,\ldots,\,E_d\;$ are known as 
the {\it primitive idempotents} of $\;\Gamma $. 
We refer to $\;E_0\;$ as the {\it trivial} idempotent.

\medskip
\noindent
Let  $\;\theta_0,\theta_1,\ldots,\theta_d$ denote the real numbers satisfying 
$
A = \sum^d_{i=0} \theta_i E_i.
$
Observe $\;AE_i\,=\,E_iA\,=\,\theta_iE_i$ for $0 \le i \le d$, and
that 
$\;\theta_0,
\,\theta_1,\,\ldots,\,\theta_d\;$ are distinct since $A$ generates $M$.
It follows from (\ref{eq:e0}) that $\;\theta_0=k$, and it is known 
$-k\le \theta_i\le k$ for $0\le i\le d$ \cite[p.197]{BI}. 
We refer to $\theta_i$ as the {\it eigenvalue} of $\Gamma $ associated 
with $E_i$, and  call $\theta_0$ the {\it trivial} eigenvalue.
For each integer $i$ $(0 \le i \le d)$, let $m_i$ denote the 
rank of $E_i$.  We  refer to $m_i$ as the {\it multiplicity} of $E_i$
(or $\theta_i$). We observe $m_0=1$.

\medskip
\noindent We now recall the cosines.
Let $\theta $ denote an eigenvalue of $\Gamma $, and  let 
$E$ denote the associated primitive idempotent. 
Let $\sigma_0, \sigma_1, \ldots, \sigma_d$ denote the real numbers satisfying
\begin{equation}
E = \vert X \vert^{-1}m\sum_{i=0}^d \sigma_iA_i,
\label{pre3}
\end{equation}
where $m$ denotes the multiplicity of $\theta $.
Taking the trace in (\ref{pre3}), we find $\sigma_0=1$.  
We often abbreviate $\sigma = \sigma_1$.
We refer to $\sigma_i$ as the 
{\it ith cosine} of $\Gamma $ with respect to $\theta $ (or $E$), and call
 $\sigma_0, \sigma_1, \ldots, \sigma_d$  
the {\it cosine sequence} of $\Gamma $ associated with $\theta $ (or $E$). 
We interpret the cosines as follows.
Let $\R^X$ denote the vector space consisting of all column vectors
with entries in $\R$ 
whose coordinates are indexed by $X$. We observe $ 
\hbox{Mat}_X(\R)$ acts on $\R^X$ by left multiplication.
We endow
$\R^X$ with the  Euclidean inner product satisfying 
\begin{equation}
\langle u,v \rangle  = u^t v \qquad \qquad (u,v \in \R^X),
\label{pre1}
\end{equation}
where $t$  denotes transposition. 
For each $x \in X$, let $\hat x $ denote the element in $\R^X$ with a 
1 in coordinate $x$, and 0 in all other coordinates.
We note
$\lbrace \hat x \;\vert \;x \in X\rbrace $ is an orthonormal basis
for $\R^X$.

\begin{lemma} \label{L4.1} \label{ZL2.1}       
Let $\;\Gamma=(X,R)\;$ denote a distance-regular graph with diameter
$\;d\geq 3$. Let $\;E\;$ denote a primitive idempotent of $\;\Gamma$,
and let $\;\sigma_0,\,\sigma_1,\,\ldots,\,\sigma_d\;$ denote
the associated cosine sequence.
Then for all integers $i$ $(0 \le i \le d)$, and for all
$\;x,\,y \in X$ such that $\partial(x,y)=i$, the following (i)--(iii)
hold.
\begin{description} \itemsep -1pt
\item{(i)} \ \ $\langle E\hat x, E\hat y\rangle =  
	    m|X|^{-1}\sigma_i $, where $m$ denotes the multiplicity of $E$. 
\item{(ii)} \ The cosine of the angle between the vectors $\;E\hat x\;$
	      and $\;E\hat y\;$ equals $\sigma_i $. 
\item{(iii)}\ 
$\;-1\le \sigma_i \leq 1 $.
\end{description} 
\end{lemma} 

\proof Line (i) is a routine application of
(\ref{eq:eimult}), (\ref{pre3}),
(\ref{pre1}). 
 Line (ii) is  immediate from (i), and (iii) is immediate from (ii).
\qed

%
%
%
%
%


%
%
%
\begin{lemma}\cite[Sect. 4.1.B]{BCN} \label{Lem4.3}  \label{ZL2.4}            
Let $\Gamma$ denote a distance-regular graph with diameter 
$\;d\geq 3$. Then for any complex numbers  
$\;\theta,\, \sigma_0,\,\sigma_1,\,\ldots,\,\sigma_d$, the following 
are equivalent.
\begin{description} \itemsep -3pt
\item{(i)} \ \ \ $\theta\;$ is an eigenvalue of $\;\Gamma $, and 
	   $\;\sigma_0,\,\sigma_1,\,\ldots,\,\sigma_d\;$ is the 
	   associated cosine sequence.
\item{(ii)} \ \ $\sigma_0=1$, and 
\begin{equation}
	     c_i\sigma_{i-1} +a_i\sigma_i+b_i\sigma_{i+1} 
	     = \theta \sigma_i \qquad \qquad (0 \le i \le d),
\label{pre2A}
\end{equation} 
             where $\;\sigma_{-1}$ and  $\sigma_{d+1}\;$ are
	     indeterminates.
\item{(iii)} \ $\sigma_0=1$, $k \sigma = \theta$, and 
            \begin{equation} 
             c_i(\sigma_{i-1} - \sigma_i) - b_i(\sigma_i-\sigma_{i+1}) 
             = k(\sigma-1) \sigma_i \qquad \qquad (1 \le i \le d),
           \label{pre2B}
	   \end{equation}
             where  $\sigma_{d+1}\;$ is an  
	     indeterminate.
	     \qed
%
\end{description}
\end{lemma}

%
%
%

\noindent 
For later use we record a number of consequences of Lemma \ref{Lem4.3}.

\begin{lemma}            \label{TECH}                  
Let $\Gamma$ denote a distance-regular graph with diameter $\;d\geq 3$.
Let $\theta\;$ denote an eigenvalue of $\;\Gamma$,  
and let $\;\sigma_0, \sigma_1,\,\ldots,\sigma_d\;$
denote the associated cosine sequence. Then (i)--(vi) hold below.
\begin{enumerate} \itemsep -5pt
\item 
$ kb_1\sigma_2 = \theta^2 -a_1\theta -k.$
\item ${\displaystyle{
kb_1(\sigma - \sigma_2) = (k-\theta)(1+\theta).
}}$
\item ${\displaystyle {kb_1(1-\sigma_2)= (k-\theta)(\theta+k-a_1)}}.$
\item $k^2b_1(\sigma^2-\sigma_2)=(k-\theta)(k+\theta (a_1+1))$.
\item $c_d(\sigma_{d-1}-\sigma_d) =k(\sigma-1)\sigma_d$.
\item $a_d(\sigma_{d-1}-\sigma_d) = k(\sigma_{d-1}-\sigma\sigma_d)$.
\end{enumerate}
\end{lemma}
\proof To get (i), 
set $i=1$ in (\ref{pre2A}), and solve for $\sigma_2$.
Lines 
(ii)--(iv) are routinely verified using (i) above and $k\sigma = \theta $.
To get  (v), set $i=d$, $b_d=0$   in Lemma \ref{ZL2.4}(iii). 
To get 
(vi), set $c_d=k-a_d$ in (v) above, and simplify the
result.
\qed

\medskip
\noindent In this article,     
the second largest
and minimal eigenvalue of a distance-regular graph turn out to be
of particular interest.      
In the next several lemmas, we give some basic information on these
eigenvalues.
\begin{lemma}\cite[Lem. 13.2.1]{God2} \label{L5.1} \label{ZL2.6}   
Let $\Gamma$ denote a distance-regular graph with diameter 
$\;d\geq 3$, and eigenvalues $\;\theta_0> \theta_1>\cdots >\theta_d$.
Let $\theta$ denote one of $\theta_1, \theta_d$ and let
$\;\sigma_0,\,\sigma_1,\,\ldots,\,\sigma_d$ denote the cosine sequence
for $\theta $.

\begin{description} \itemsep -1pt
\item{(i)} \ \ Suppose $\theta=\theta_1$. Then $\;\sigma_0 \;>\;\sigma_1\;>\;\cdots\;>\;\sigma_d$.
\item{(ii)} \ Suppose $\theta=\theta_d$. Then  $\;\; (-1)^i\sigma_i > 0$ \qquad \qquad $(0 \le i \le d)$. \qed
\end{description}
\end{lemma} 

\medskip \noindent 
Recall a distance-regular graph  $\Gamma $ is {\it bipartite}
whenever the intersection numbers satisfy $\;a_i=0\;$ for 
$0 \le i \le d$, where $d$ denotes the diameter. 

\begin{lemma}                                          \label{drgbip}
Let $\;\Gamma=(X,R)\;$ denote a  distance-regular graph
with diameter $\;d\geq 3$. Let $\theta_d$ denote the minimal eigenvalue 
of $\Gamma $, and let $\sigma_0, \sigma_1,\ldots,\sigma_d$ denote the 
associated  cosine sequence.  Then the following are equivalent:
%
(i) $\Gamma $ is bipartite; 
(ii) $\theta_d=-k$;
(iii) $\sigma_1=-1$; 
(iv) $\sigma_2=1$. 
Moreover, suppose (i)--(iv) hold. Then 
$
\;\sigma_i = (-1)^i \;$ for $\; 0 \le i \le d.$
\end{lemma}

\proof The equivalence of  (i), (ii) follows from 
\cite[Prop. 3.2.3]{BCN}.
The equivalence of (ii), (iii) is immediate from $k \sigma_1=\theta_d$.
The remaining implications follow from \cite[Prop. 4.4.7]{BCN}.  \qed

\begin{lemma} \label{ZL2.0}       
Let $\Gamma$ denote a distance-regular graph with
 diameter $\;d\geq 3$ and eigenvalues 
$\;\theta_0 > \theta_1 > \cdots > \theta_d$.  Then  (i)--(iii)
hold below.
\begin{description} \itemsep -1pt
\item{(i)}   \ \ $0 < \theta_1 < k$.
\item{(ii)} \    $a_1-k\le \theta_d <-1$.
\item{(iii)}\ Suppose $\Gamma $ is not bipartite. 
              Then $\,a_1 -k< \theta_d$. 
\end{description} 
\end{lemma}

\proof (i) The eigenvalue $\theta_1$ is positive by 
           \cite[Cor. 3.5.4]{BCN}, and we have seen $\theta_1 < k$.  


\noindent (ii)
 Let $\sigma_1, \sigma_2$ denote
the first and second cosines for $\theta_d$. Then $\sigma_2 \le 1$
by Lemma \ref{ZL2.1}(iii), so $a_1-k\leq \theta_d$ in view of 
Lemma \ref{TECH}(iii).  Also
$\sigma_1 < \sigma_2$ by Lemma \ref{ZL2.6}(ii), so $\theta_d <-1$ in view of 
Lemma \ref{TECH}(ii).

\noindent (iii) 
Suppose $\theta_d=a_1-k$. Applying
 Lemma \ref{TECH}(iii), we find
$\sigma_2=1$, where $\sigma_2$ denotes the second cosine for $\theta_d$. 
Now  $\Gamma $ is bipartite
by Lemma \ref{drgbip}, contradicting our assumptions.  Hence 
$\theta_d > a_1-k$, as desired.
\qed

%
%
%
%
%
%
%
%

\begin{lemma} \label{AL2.3}                               
Let $\;\Gamma=(X,R)\;$ denote a nonbipartite distance-regular graph
with diameter $\;d\geq 3$, let $x,y$ denote adjacent vertices in 
$X$, and let $\;E\;$ denote a nontrivial primitive 
idempotent of $\;\Gamma$. Then the vectors $\;E\hat x$ and $\;E\hat y$
are linearly independent. 
\end{lemma} 

\proof
Let $\sigma $ denote the first cosine associated to $E$. Then 
$\sigma \not= 1$, since $E$ is nontrivial, and $\sigma \not= -1$, 
since $\Gamma$ is not bipartite. Applying Lemma \ref{ZL2.1}(ii), we see
 $\;E\hat x$ and $\;E\hat y$ are linearly independent.  \qed

%
%
%
%
%
%
%
%
%



\medskip
\noindent We mention a few results on the intersection numbers.

\begin{lemma} \cite[Prop. 5.5.1]{BCN}  \label{ZZL2.11}          
Let $\Gamma$ denote a distance-regular graph with diameter $d\ge 3$
and $a_1 \not= 0$. Then
$
\;a_i \not= 0  \;(1 \le i \le d-1)$.
\qed
\end{lemma}


%

%
\begin{lemma}\cite[Lem. 4.1.7]{BCN} \label{AL2.13}                             
Let $\Gamma$ denote a distance-regular graph with 
 diameter $\;d\geq 3$. Then the intersection numbers satisfy  
$$
 p_{ii}^1 =
                            {b_1b_2\dots b_{i-1} \over c_1c_2\dots c_i}a_i,
\qquad \qquad 
 p_{i-1,i}^1 =  
                           {b_1b_2\dots b_{i-1} \over c_1c_2\dots c_{i-1}}
			    \qquad \qquad (1 \leq i \leq d). 
$$
%
%
\qed
\end{lemma}



\medskip
\noindent For the remainder of this section, we describe a point of
view we will adopt throughout the paper.

\begin{df} \label{ZD4.1}                            
Let $\Gamma = (X, R)$ denote a distance-regular graph with diameter
$d\geq 3$, and fix adjacent 
vertices $x, y \in X$.  For all integers $i$ and $j$ we define 
$D_i^j = D_{i}^{j}(x,y)$ by
\begin{equation}
D_i^j \;=\; \Gamma_i(x) \cap \Gamma_j(y).   
                                            \qquad \label{ZD4.1A}
\end{equation}
We observe $|D_i^j|=p_{ij}^1\;$ for $0\leq i,j\leq d$, and $D_i^j=\emptyset $
otherwise. We visualize the $D_i^j$ as follows.
\end{df}

\bigskip
\input psfig.sty
\centerline{\psfig{figure=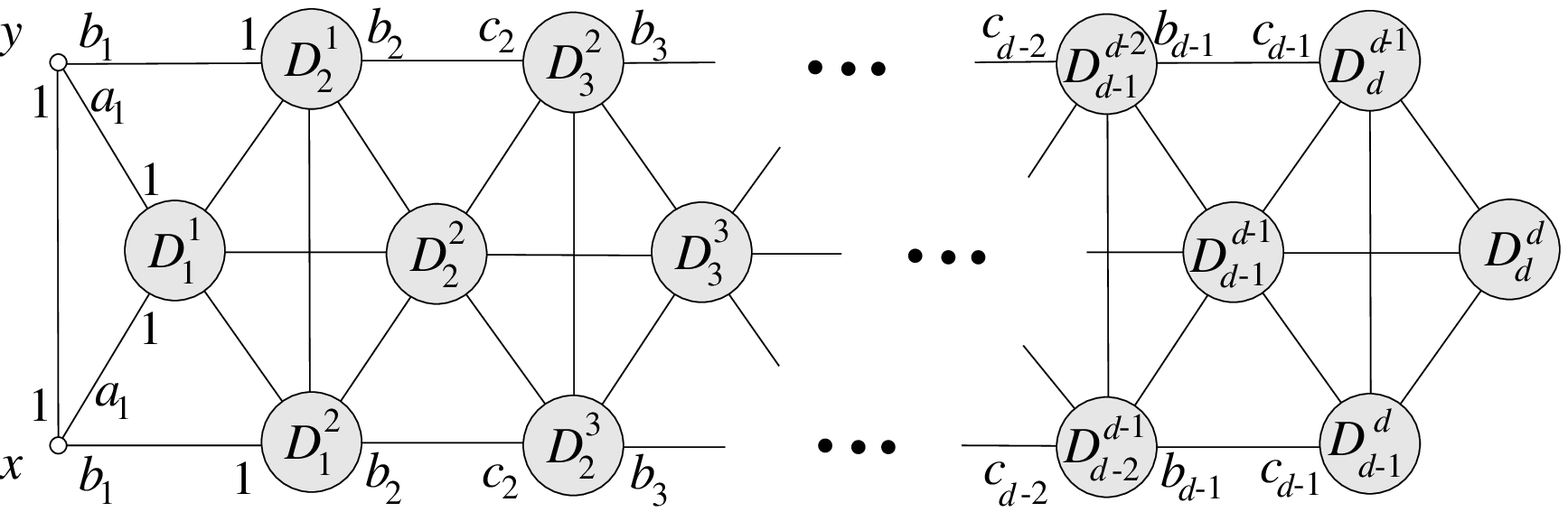,height=4.5cm}} 
\vskip .5cm 
\baselineskip 7pt \par{\leftskip 1cm \rightskip 1cm \noindent 
\small Figure 2.1: Distance distribution corresponding to an edge.
Observe: $D_i^{i-1} \cup D_i^i \cup D_i^{i+1} = \Gamma_i(x)$ for
$i=1,\dots,d$. The number beside edges connecting cells $D_i^j$ indicate
how many neighbours a vertex from the closer cell has in the other cell,
see Lemma \ref{NNL6.1}.

\par} \baselineskip=\normalbaselineskip

\begin{lemma} \label{NL4.6}  \label{NNL6.1}             
Let $\;\Gamma=(X,R)\;$ denote a distance-regular graph with diameter 
$\;d\geq 3$. Fix adjacent vertices $x, y \in X$, and 
pick any integer $i$ $ (1\leq i \leq d)$. Then with reference to 
Definition \ref{ZD4.1},
 the following (i) and (ii)
hold.
\begin{description} \itemsep -3pt
\item{(i)} Each $z\in  D_{i-1}^{i}$
           (resp. $D_{i}^{i-1}$)
           is adjacent to 
           \begin{description} \itemsep -3pt
           \item{(a)} precisely 
		      \ \ \ \ \ \ \ \ \ \ \ \ \ \ \ \ \ \ 
		      $c_{i-1}$ 
		      \ \ \ \ \  \ \ \ \ \ \ \ \ \ \
		      vertices in $D_{i-2}^{i-1}$
                      (resp. $D_{i-1}^{i-2}$),
           \item{(b)} precisely 
		      \ \ \ \ \ \ 
		      $c_i-c_{i-1} -|\Gamma(z)\cap D_{i-1}^{i-1}|$
		      \ \ 
                      vertices in $D_{i}^{i-1}$
                      (resp. $D_{i-1}^{i}$),
           \item{(c)} precisely 
		      \ \ \ \ \ \ \ \ \ 
		      $a_{i-1} -|\Gamma(z)\cap D_{i-1}^{i-1}|$
		      \ \ \ \ 
                      vertices in $D_{i-1}^{i}$
                      (resp. $D_{i}^{i-1}$),
           \item{(d)} precisely 
		      \ \ \ \ \ \ \ \ \ \ \ \ \ \ \ \ \ \ \ 
		      $b_i$ 
		      \ \ \ \ \ \ \ \ \ \ \ \ \ \ \ \ 
		      vertices in $D_{i}^{i+1}$
                      (resp. $D_{i+1}^{i}$),
           \item{(e)} precisely 
		      \ \ \ \ \ \ 
		      $a_i-a_{i-1} +|\Gamma(z)\cap D_{i-1}^{i-1}|$
		      \ \ 
                      vertices in $D_{i}^{i}$.
           \end{description}
\item{(ii)} Each $z \in D_{i}^{i}$ 
           is adjacent to 
           \begin{description} \itemsep -3pt
           \item{(a)} precisely 
		      \ \ \ \ \ \ \ \ \ \ \ \ \ \ \ \ \ 
		      $c_i - |\Gamma(z) \cap D_{i-1}^{i-1}|$
		      \ \ \ \ \ \ \ \ \ \ \ \ \ \ \ \ \ \ \ 
                      vertices in $D_{i-1}^{i}$,
           \item{(b)} precisely 
		      \ \ \ \ \ \ \ \ \ \ \ \ \ \ \ \ \ 
		      $c_i - |\Gamma(z) \cap D_{i-1}^{i-1}|$
		      \ \ \ \ \ \ \ \ \ \ \ \ \ \ \ \ \ \ \ 
                      vertices in $D_{i}^{i-1}$,
           \item{(c)} precisely 
		      \ \ \ \ \ \ \ \ \ \ \ \ \ \ \ \ \ 
		      $b_i - |\Gamma(z) \cap D_{i+1}^{i+1}|$
		      \ \ \ \ \ \ \ \ \ \ \ \ \ \ \ \ \ \ \ 
                      vertices in $D_{i}^{i+1}$,
           \item{(d)} precisely 
		      \ \ \ \ \ \ \ \ \ \ \ \ \ \ \ \ \ 
		      $b_i - |\Gamma(z) \cap D_{i+1}^{i+1}|$
		      \ \ \ \ \ \ \ \ \ \ \ \ \ \ \ \ \ \ \ 
                      vertices in $D_{i+1}^{i}$,
           \item{(e)} precisely 
		      \ \ \ \ \ \ 
		      $a_i-b_i-c_i + |\Gamma(z) \cap D_{i-1}^{i-1}|
                      + |\Gamma(z) \cap D_{i+1}^{i+1}|$
                      vertices in $D_{i}^{i}$.
           \end{description}
\end{description}
\end{lemma}

\proof Routine. \qed

\section{Edges that are tight with respect to an eigenvalue}     

Let $\Gamma =(X,R)$ denote a graph, and let $\Omega $ denote
a nonempty subset of $X$. By the {\it vertex subgraph} of 
$\Gamma $ {\it induced on } $\Omega $, we mean the graph with
vertex set $\Omega $, 
and edge set $\lbrace xy \;\vert \;x, y \in \Omega,\quad 
xy \in R\rbrace $.

\begin{df} \label{ZL3.0A}                              
Let $\Gamma = (X, R)$ denote a distance-regular graph with diameter
$d \geq 3$ and intersection number $a_1 \neq 0$. 
For each edge $xy \in R$, we define the scalar $f=f(x,y)$ by
\begin{equation}                                     \label{fdef} 
f :=  a_1^{-1}\Bigl\vert\{ (z,w) \in X^2 \mbox{ $|$ } 
z,w \in D^1_1,\mbox{ } \dist(z, w) = 2 \}\Bigr\vert,
\end{equation}
where $D^1_1=D^1_1(x,y)$ is from (\ref{ZD4.1A}). 
We observe $f$ is the average valency  of the complement of the 
vertex subgraph induced  on $D^1_1$.
\end{df}

\noindent 
We begin with some elementary facts about $f$.

\begin{lemma} \label{ZL3.0}                              
Let $\Gamma = (X, R)$ denote a distance-regular graph with diameter
$d \ge 3$ and $a_1 \ne 0$. Let  $x,y$ denote adjacent vertices in 
$X$. Then with reference to (\ref{ZD4.1A}), (\ref{fdef}), 
lines (i)--(iv) hold below.
\begin{description} \itemsep -3pt
\item{(i)} \ \  The number of edges in $R$ connecting a vertex
	       in $D_1^1$ with a vertex in 
	        $D_1^2$
		is equal to  $a_1f$.
\item{(ii)} \   The number of edges in the vertex subgraph induced on 
	        $D^1_1$ is equal to $a_1(a_1-1-f)/2$.
\item{(iii)}    The number of edges in the vertex subgraph induced 
	        on $D_1^2$ is equal to $a_1(b_1-f)/2$.
\item{(iv)}  \ $0\le f, \qquad  f \le a_1-1, \qquad  f \le b_1.$
\end{description}
\end{lemma}

\proof Routine.  \qed

\medskip\noindent
The following lemma provides another bound for $f$.

\begin{lemma} \label{ZL3.1}                              
Let $\Gamma = (X,R)$ denote a distance-regular graph with
diameter $d \geq 3$ and $a_1 \neq 0$. Let $x,y$ denote adjacent
vertices in $X$, 
 and write $f=f(x,y)$. Then  for each 
nontrivial eigenvalue
$\theta$ of $\Gamma$, 
\begin{eqnarray}
(k+\theta)(1+\theta)\;f \;\le\;  b_1 \big( k + \theta(a_1+1)\big).
                                                   \label{ZL3.1A}
\end{eqnarray}
\end{lemma}

\proof
Let $\sigma_0,\dots,\sigma_d$ denote the cosine sequence of
$\theta$ and let $E$ denote the corresponding primitive idempotent.
Set 
$$
w := \sum_{z \in D_1^1}  \hat{z}, 
$$
where $D^1_1=D^1_1(x,y)$ is from (\ref{ZD4.1A}). 
Let $G$ denote the Gram matrix for the
vectors $E \hat{x}$, $E \hat{y}$, $Ew$; that is 
$$
G:= \left( \begin{array}{ccc} 
\| E \hat{x} \|^2                                 & 
\mbox{$\langle E \hat{x},$ $E \hat{y} \rangle$} &
\mbox{$\langle E \hat{x},$ $Ew \rangle$}           \cr
\mbox{$\langle E \hat{y},$ $E \hat{x} \rangle$} &
\| E \hat{y} \|^2                                 & 
\mbox{$\langle E \hat{y},$ $Ew \rangle$}           \cr
\mbox{$\langle Ew,$ $E \hat{x} \rangle$}           &
\mbox{$\langle Ew,$ $E \hat{y} \rangle$}           &
\| Ew \|^2                                           \cr
 \end{array}\right).  
$$
On one hand, the matrix $G$ is positive semi-definite, so it has 
nonnegative determinant. On the other hand, by  Lemma \ref{ZL2.1}, 
\begin{eqnarray*}
{\rm det}(G) \;&=&\; {m^3 |X|^{-3}} \;
                    \hbox{det} \left( \begin{array}{ccc} 
                    \sigma_0 & \sigma_1 & a_1 \sigma_1 \cr
                    \sigma_1 & \sigma_0 & a_1 \sigma_1 \cr
                    a_1 \sigma_1 & a_1 \sigma_1 & 
         a_1 \big(\sigma_0 + (a_1- f -1) \sigma_1 + f \sigma_2\big) \cr
                    \end{array}\right)  \\
             \;&=&\; {m^3 a_1 |X|^{-3}}\; (\sigma-1) \;
                    \Big((\sigma-\sigma_2)(1+\sigma)f
                         \,-\,(1-\sigma)(a_1\sigma+1+\sigma)\Big),
\end{eqnarray*}
where $m$ denotes the multiplicity of $\theta$.
Since $a_1>0$ and $\sigma < 1$, we find 
\begin{equation}                                         \label{ZL3.1B}
(\sigma-\sigma_2)(1+\sigma)\; f \;\le\; (1-\sigma)(a_1\sigma+1+\sigma).
\end{equation}
Eliminating $\sigma,\,\sigma_2$ in (\ref{ZL3.1B}) using $\theta=k\sigma$
and Lemma \ref{TECH}(ii), and simplifying the result using $\theta < k$,
we routinely obtain (\ref{ZL3.1A}). \qed

\begin{cor}\label{lembf3} \label{NL3.4}\label{ZC3.2} 
Let $\Gamma = (X,R)$ denote a distance-regular graph with diameter 
$d \ge 3$ and $a_1 \ne 0$. Let $x,y$ denote adjacent vertices in $X$, 
and let $\theta$ denote a nontrivial eigenvalue of $\Gamma$. Then with 
reference to Definition \ref{ZD4.1}, the following are equivalent.
\begin{description}\itemsep 30pt \parskip -30pt 
\item{(i)} \ \ \  Equality is attained in (\ref{ZL3.1A}).
\item{(ii)} \ \   $E \hat{x},$ $E \hat{y},$ ${\displaystyle
                  \sum_{z \in D_1^1} E \hat{z} }$
	          are linearly dependent.
\item{(iii)} \ $\displaystyle \sum_{z \in D_1^1} E \hat{z}
                = {a_1 \theta \over k+\theta} (E \hat{x} + E \hat{y})$.
\end{description}
\end{cor}

\medskip\noindent
We say the edge $xy$ is {\bf tight with respect to $\theta$} whenever
(i)--(iii) hold above.

\proof (i)\iff(ii) Let the matrix $G$ be as in the proof of Lemma  
\ref{ZL3.1}.  Then we find (i) holds if and only if $G$ is singular,  
if and only if (ii) holds.

\noindent 
(ii) \implies (iii) $\Gamma$ is not bipartite since $a_1 \not= 0$,
so $E \hat{x},$ and $E \hat{y}$ are linearly independent by 
Lemma \ref{AL2.3}. It follows 
\begin{equation}
    \sum_{z \in D_1^1} E \hat{z}
\; = \;  \alpha E \hat x + \beta E \hat y
\label{Sdependence}
\end{equation}
for some $\,\alpha, \beta \in \R$.
Taking the inner product of (\ref{Sdependence})  with each of  
$E \hat x$,  $E \hat y$ using Lemma \ref{ZL2.1},
  we readily obtain 
$\alpha =\beta = a_1\theta (k+\theta)^{-1}$.

\noindent (iii) \implies (ii) Clear. \qed

\bigskip\noindent
Let $\Gamma=(X,R)$ denote a distance-regular graph with diameter 
$d\ge 3$, $a_1 \ne 0$, and eigenvalues 
$\theta_0 >\theta_1 >\cdots >\theta_d$.
Pick adjacent vertices $x,y \in X$, and write $f=f(x,y)$. 
Referring to (\ref{ZL3.1A}), we now consider which of 
$\theta_1,\theta_2,\ldots, \theta_d$ gives the best  bounds for $f$.
Let $\theta $ denote one of $\theta_1,\theta_2,\ldots, \theta_d$. 
Assume $\theta \not=-1$; otherwise (\ref{ZL3.1A}) gives no information 
about $f$. If $\theta > -1$ (resp. $\theta < -1$), line (\ref{ZL3.1A}) 
gives an upper (resp. lower) bound for $f$.
Consider the partial fraction  decompostion
$$
b_1 {k + \theta(a_1 + 1) \over (k+\theta)(1+\theta)} \ \ =\ \
{b_1 \over k-1} \Big({ka_1\over k+\theta}+{b_1 \over 1+\theta}\Big).
$$
Since the map $F: \rr \setminus \{-k,-1\} \longrightarrow \rr$, defined
by
$$
x \mapsto {ka_1 \over k+x} + {b_1 \over 1 +x}
$$ 
is strictly decreasing on the intervals $(- k, -1)$ and $(-1, \infty)$,
we find in view of Lemma \ref{ZL2.0} that
the least  upper bound for $f$ is obtained at $\theta=\theta_1$, and   
the greatest  lower bound is obtained at $\theta=\theta_d$.

\begin{th}\label{thmbf} \label{NT3.9} \label{ZT3.3} 
Let $\Gamma = (X,R)$ denote a distance-regular graph with diameter 
$d \ge 3$,  $a_1 \ne 0$, and eigenvalues 
$\theta_0 >\theta_1 >\cdots >\theta_d $. For all edges $xy \in R$,
\begin{equation}                                         \label{ZT3.3A}
b_1 {k+ \theta_d(a_1+1) \over (k+\theta_d)(1+\theta_d)} 
\ \ \le \ \ f(x,y) \ \ \leq \ \ 
b_1 {k+ \theta_1(a_1+1) \over (k+\theta_1)(1+\theta_1)}.
\end{equation}
\end{th}

\proof This is immediate from (\ref{ZL3.1A}) and Lemma \ref{ZL2.0}.\qed

\begin{cor}\label{thmbf2} \label{NT3.10}     
Let $\Gamma = (X,R)$ denote a distance-regular graph with diameter
$d \ge 3$, $a_1 \ne 0$, and  eigenvalues
$\; \theta_0 > \theta_1 > \cdots > \theta_d$. For all edges $xy \in R$,
\begin{description} \itemsep -3pt
\item{(i)} \ \ \ $xy$ is tight with respect to $\theta_1$ if and only if
	     equality holds in the right inequality of (\ref{ZT3.3A}), 
\item{(ii)} \ \  $xy$ is tight with respect to $\theta_d$ if and only if
	     equality holds in the left inequality of
            (\ref{ZT3.3A}), 
\item{(iii)} \ $xy$ is not tight with respect to $\theta_i$ for 
             $2 \le i \le d-1$. 
\end{description}
\end{cor}
  
\proof
(i),(ii) Immediate from  (\ref{ZL3.1A}) and Corollary \ref{ZC3.2}.

\noindent 
(iii) First suppose $\theta_i=-1$. We do not have equality for 
$\theta=\theta_i$ in (\ref{ZL3.1A}), since the left side  equals $0$, 
and the right side equals $b_1^2$.  
In particular, $xy$ is not tight with respect to $\theta_i$. 
Next suppose $\theta_i\not=-1$. Then we do not have equality for
$\theta=\theta_i$ in (\ref{ZL3.1A}) in view of the above mentioned fact,
that the function $F$ is strictly decreasing on the intervals 
$(- k,-1)$ and  $(-1,\infty)$. \qed

\section{Tight edges and combinatorial regularity}     

\begin{th} \label{ZZT4.2}                                
Let $\;\Gamma=(X,R)\;$ denote a distance-regular graph with diameter 
$\;d\geq 3$ and intersection number $a_1\not=0$.
Let $\theta$ denote a nontrivial eigenvalue of
$\Gamma$, and let $\sigma_0,\sigma_1,\dots,\sigma_d$ denote its cosine
sequence. 
Let $x,y $ denote adjacent vertices in $X$. Then with reference to 
Definition \ref{ZD4.1}, 
the following are equivalent.
\begin{description} \itemsep -3pt
\item{(i)} \ \ \ $xy$ is tight with respect to $\theta$.
\item{(ii)} \ \ For $1 \le i \le d$; both  
                $\sigma_{i-1} \not= \sigma_i$, and for all 
                $z\in D_{i-1}^{i}$
            \begin{eqnarray}
            |\Gamma_{i-1}(z) \cap D_1^1| \;&=&\; {a_1 \over 1+\sigma}\; 
	    {\displaystyle
            {\sigma\sigma_{i-1}-\sigma_i \over \sigma_{i-1}-\sigma_i}},
				   	       \label{ZZT4.2A}
            \\
            |\Gamma_i(z) \cap D_1^1| \;&=&\; {a_1 \over 1+\sigma}\;
	    {\displaystyle
            {\sigma_{i-1}-\sigma\sigma_i \over \sigma_{i-1}-\sigma_i}}.
						   \label{ZZT4.2B}
            \end{eqnarray}
\end{description}
\end{th}

\proof
(i) \implies (ii) Let the integer $i$ be given. Observe
 by Corollary \ref{NT3.10} that 
$\theta $ is either 
the second largest eigenvalue $\theta_1$ or the least eigenvalue
 $\theta_d$,  
 so 
$\sigma_{i-1}\not=\sigma_i$ in view of  Lemma \ref{ZL2.6}. Pick any
$z\in D_{i-1}^{i}$.
Observe 
 $D_1^1$ contains $a_1$ vertices, and each is at distance 
 $i-1$ or $i$ from $z$, so
%
\begin{equation}
|\Gamma_{i-1}(z) \cap D_1^1| + |\Gamma_{i}(z) \cap D_1^1| \;=\; a_1.
                                               \label{ZZT4.2C}
\end{equation}
Let $E$ denote the primitive idempotent associated to $\theta$.
By Corollary \ref{ZC3.2}(iii), and since $xy$ is tight with respect to 
$\theta$,
\begin{equation}
\sum_{w \in D_1^1} E \hat{w}
\;=\; {a_1 \sigma \over 1+\sigma} (E \hat{x} + E \hat{y}).
\label{Innprod}
\end{equation}
Taking the inner product of
(\ref{Innprod}) 
with 
$E \hat z$ using 
Lemma \ref{ZL2.1},
 we obtain
\begin{equation}
\sigma_{i-1} |\Gamma_{i-1}(z) \cap D_1^1| +
\sigma_i |\Gamma_{i}(z) \cap D_1^1| \;=\;
{a_1 \sigma \over 1 + \sigma} (\sigma_{i-1} + \sigma_i).
                                               \label{ZZT4.2D}
\end{equation}
%
Solving the system (\ref{ZZT4.2C}), (\ref{ZZT4.2D}), we routinely obtain
(\ref{ZZT4.2A}), (\ref{ZZT4.2B}).
\\
(ii) \implies (i)
We show equality holds in (\ref{ZL3.1A}). 
Counting the edges between $D_1^1$ and $D_1^2$ using (\ref{ZZT4.2A}) (with
$i=2$), we find in view of Lemma \ref{ZL3.0}(i) that 
\begin{equation} 
f(x,y) \;=\; b_1 {\sigma^2 - \sigma_2 \over (1+\sigma)(\sigma-\sigma_2)}.
\label{eq:fvalue}
\end{equation}
Eliminating $\sigma, \sigma_2$ in (\ref{eq:fvalue}) using 
$\theta=k \sigma $ and Lemma \ref{TECH}(ii),(iv),  
 we  readily find
 equality holds in 
(\ref{ZL3.1A}). Now $xy$ is tight with respect to $\theta$ by Corollary
\ref{ZC3.2}. \qed

\begin{th} \label{NT4.12} \label{UZT4.3}                
Let $\;\Gamma=(X,R)\;$ denote a distance-regular graph with diameter 
$\;d\geq 3$ and $a_1\not=0$.  Let $\theta$ denote a nontrivial eigenvalue of
$\Gamma$, and let $\sigma_0,\sigma_1,\dots,\sigma_d$ denote its cosine
sequence.
Let $x,y $ denote adjacent vertices in $X$. Then with reference to 
Definition \ref{ZD4.1}, 
 the following are equivalent.
\begin{description} \itemsep -3pt
\item{(i)} \ \ \ $xy$ is tight with respect to $\theta$,

\item{(ii)} \ For  $1 \le i \le d-1$; both 
            $\sigma_i \not= \sigma_{i+1}$, and for all $z\in D_i^i$
            \begin{eqnarray}
      |\Gamma_{i+1}(z) \cap D_1^1| \;&=&\; |\Gamma_{i-1}(z) \cap D_1^1|\;
           {\displaystyle 
	   {\sigma_{i-1} - \sigma_i \over \sigma_i - \sigma_{i+1}}  
           \;+\; a_1\; {1-\sigma \over 1+\sigma} \;
           {\sigma_i \over \sigma_i - \sigma_{i+1}}},
                                              \label{NT4.12A}
\\[1ex]
         |\Gamma_i(z) \cap D_1^1| \;&=&\; -|\Gamma_{i-1}(z) \cap D_1^1|\;
           {\displaystyle 
           {\sigma_{i-1} - \sigma_{i+1} \over \sigma_i - \sigma_{i+1}}
           \;+\; a_1\; {2\sigma \over 1+\sigma}}  \nonumber \\
           &&\;-\; {\displaystyle a_1\; {1-\sigma \over 1+\sigma} \;
             {\sigma_{i+1}\over \sigma_i - \sigma_{i+1}}}.
                                                  \label{NT4.12B}
           \end{eqnarray}
\end{description}
Suppose (i)--(ii) above, and that $a_d\not= 0$. 
Then for all $z\in D_d^d$
\begin{eqnarray}
|\Gamma_{d-1}(z) \cap D_1^1| \;&=& \;-\; a_1 {1-\sigma \over 1+\sigma} \;
                                {\sigma_d\over \sigma_{d-1} - \sigma_d},
                                                  \label{NT4.12C} 
\\
|\Gamma_{d}(z)\cap D_1^1| \;&=&\; a_1\;+\; a_1{1-\sigma\over 1+\sigma} \;
                                {\sigma_d\over \sigma_{d-1} - \sigma_d}.
                                                  \label{NT4.12D}
\end{eqnarray}
\end{th}

\proof (i) \implies (ii)
 Let the integer $i$ be given. Observe 
by Corollary \ref{NT3.10} that 
 $\theta $ is either the second largest eigenvalue 
 $\theta_1$ or the least eigenvalue  
$\theta_d$,  
 so 
$\sigma_{i}\not=\sigma_{i+1}$ by Lemma \ref{ZL2.6}. Pick any $z\in D_i^{i}$.
Proceeding as in the proof of Theorem \ref{ZZT4.2} (i) \implies (ii),
we find 
\begin{eqnarray}
|\Gamma_{i-1}(z) \cap D_1^1| +|\Gamma_{i}(z) \cap D_1^1| +
|\Gamma_{i+1}(z) \cap D_1^1|\;&=&\; a_1,
						   \label{NT4.12F}
\\
\sigma_{i-1} |\Gamma_{i-1}(z) \cap D_1^1| + 
\sigma_i |\Gamma_{i}(z) \cap D_1^1| + \sigma_{i+1} 
|\Gamma_{i+1}(z)\cap D_1^1| \;&=&\; {2\sigma \sigma_i  a_1 \over 1+\sigma}.
						   \label{NT4.12G}
\end{eqnarray}
Solving (\ref{NT4.12F}), (\ref{NT4.12G}) for 
$|\Gamma_i(z) \cap D_1^1|$, $|\Gamma_{i+1}(z) \cap D_1^1|$, 
we routinely  obtain (\ref{NT4.12A}) and (\ref{NT4.12B}).
\\
(ii) \implies (i) 
Setting $i=1$ in (\ref{NT4.12A}), and evaluating the result using (\ref{fdef}),
we find
\begin{equation}
f(x,y)\;=\;  
           {\displaystyle 
	   {1 - \sigma \over \sigma - \sigma_2}  
           \;+\; a_1\; {1-\sigma \over 1+\sigma} \;
           {\sigma \over \sigma - \sigma_2}}.
                                              \label{eq:fvalue2}
\end{equation}
Eliminating $\sigma, \sigma_2$ in (\ref{eq:fvalue2}) using 
$\theta = k \sigma $ and 
 Lemma \ref{TECH}(ii), we find 
 equality holds in 
(\ref{ZL3.1A}).   Now $xy$ is tight with respect to $\theta$ by
Corollary \ref{ZC3.2}.
\\
Now suppose (i)--(ii) hold above, and that $a_d \not= 0$. Pick any 
$z\in D^d_d$. Proceeding as in the proof of Theorem \ref{ZZT4.2}
(i) \implies (ii), we find 
\begin{eqnarray}
|\Gamma_{d-1}(z) \cap D_1^1| +|\Gamma_{d}(z) \cap D_1^1| \;&=&\; a_1,
						   \label{NT4.12H}
\\
\sigma_{d-1} |\Gamma_{d-1}(z) \cap D_1^1| + 
\sigma_d |\Gamma_{d}(z) \cap D_1^1|
\;&=&\; {2\sigma_d \sigma a_1 \over 1+\sigma}.
						   \label{NT4.12I}
\end{eqnarray}
Observe
$\sigma_{d-1}\not= \sigma_d$ by (ii) above,  
so the linear system
(\ref{NT4.12H}), (\ref{NT4.12I}) has unique solution 
(\ref{NT4.12C}), (\ref{NT4.12D}).\qed

\section{The tightness of an edge}       

\begin{df}\label{deftight}                       
Let $\Gamma = (X, R)$ denote a distance-regular graph with diameter
$d \geq 3$, intersection number  $a_1 \not= 0$, and eigenvalues
 $\theta_0>\theta_1>\cdots > \theta_d$.
For each edge $xy \in R$, let $t = t(x,y)$
denote the number of nontrivial eigenvalues of $\Gamma$ with respect to
which $xy$ is tight. We call $\,t$ the {\bf tightness} of the edge $xy$.
In view of Corollary \ref{NT3.10} we have:
\begin{description} \itemsep -3pt
\item{(i)}\ \ \ $t=2$ if $xy$ is tight with respect to both $\theta_1$ 
	       and $\theta_d$;          
\item{(ii)}\ \  $t=1$ if $xy$ is tight with respect to exactly one of 
                $\theta_1$ and $\theta_d$; 
\item{(iii)}\ $t=0$ if $xy$ is not tight with respect to $\theta_1$ or
$\theta_d$. 
%
\end{description}
\end{df}

\begin{th} \label{NL4.2}                        
Let $\;\Gamma=(X,R)\;$ denote a distance-regular graph with diameter 
$\;d\geq 3$ and  $a_1 \not= 0$. For all edges $xy \in R$, the
tightness 
$ t\:=\:t(x,y)$ is given by  
\begin{equation}
t= 3d+1 - \dim(MH),        \qquad \label{NL4.2A}
\end{equation}
where $M$ denotes the Bose-Mesner algebra of $\Gamma $, where 
\begin{equation}
H\;=\; \mbox{\Span} \biggl\lbrace\hat x,\ \hat y, \;
\sum_{ z \in D_1^1(x,y)} \hat z\biggr\rbrace,
 	                                  \qquad \label{NL4.2B}
\end{equation}
and where $MH$ means  
$\mbox{\Span} \lbrace mh\;|\;m\in M,\;\;h\in H\rbrace$.
%
%
%
%
\end{th}
\proof
Since $E_0, E_1, \ldots, E_d$ is a basis for $M$, and in view
of (\ref{eq:eimult}),  
$$
MH = \sum_{i=0}^d E_i H \ \ \ \ \ \ \ \ \ \ {\rm (direct \  sum),}
$$
and it follows
$$
\dim MH = \sum_{i=0}^d \dim E_i H.
$$
Note that $\dim E_0 H=1$. For $1 \leq i \leq d$, we find 
by Lemma \ref{AL2.3} and Corollary \ref{ZC3.2}(ii) that 
 $\dim E_i H=2$ if $xy$ is tight with respect to 
$\theta_i$, and $\dim E_i H=3$ otherwise. The result follows. \qed

\section{Tight graphs and the Fundamental Bound}    

In this section, we obtain an inequality involving the second
largest and minimal eigenvalue of a distance-regular graph.
To obtain it, we need the following lemma.
\begin{lemma}\label{ZT3.50}               
Let $\Gamma$ denote a nonbipartite distance-regular graph with diameter 
$d \geq 3$, and  eigenvalues 
$\;\theta_0 > \theta_1 > \cdots > \theta_d$.
Then
\begin{eqnarray}
{k+ \theta_1(a_1+1) \over (k+\theta_1)(1+\theta_1)} 
&-&
{k+ \theta_d(a_1+1) \over (k+\theta_d)(1+\theta_d)} 
                                                        \label{ZT3.5C}
\\
&& \qquad =\;\;
\Psi \,{{(a_1+1)(\theta_d-\theta_1) }\over
{(1 + \theta_1)(1+\theta_d)(k+\theta_1)(k + \theta_d)}},  \qquad \qquad  
 \label{Deltacor}
\end{eqnarray}
where
\begin{eqnarray}
 \Psi &=&\left(\theta_1 + {k \over a_1+1}\right)\left(\theta_d + {k \over a_1+1}
\right)
\;+\; {ka_1b_1 \over (a_1+1)^2}.                       \label{ZT3.5D}
\end{eqnarray}
\end{lemma}

\proof Put (\ref{ZT3.5C}) over a common denominator, and simplify. 
\qed

\smallskip
\noindent We  now present our inequality. We give two
versions.

\begin{th}\label{ZT3.5}               

Let $\Gamma$ denote a  distance-regular graph with diameter 
$d \geq 3$, and  eigenvalues 
$\;\theta_0 > \theta_1 > \cdots > \theta_d$.
Then (i), (ii) hold below.
\begin{description} \itemsep -3pt
\item{(i)} Suppose $\Gamma$ is not bipartite. Then
\begin{eqnarray}
{k+ \theta_d(a_1+1) \over (k+\theta_d)(1+\theta_d)} 
\leq 
{k+ \theta_1(a_1+1) \over (k+\theta_1)(1+\theta_1)}. 
                                                        \label{ZT3.5B}
\end{eqnarray}
\item{(ii)}
\begin{eqnarray}
 \left(\theta_1 + {k \over a_1+1}\right)\left(\theta_d + {k \over a_1+1}
\right)
\ge - {ka_1b_1 \over (a_1+1)^2}.                       \label{ZT3.5A}
\end{eqnarray}
\end{description}
We refer to (\ref{ZT3.5A}) as the {\bf Fundamental Bound}.
\end{th}

\proof (i) First assume $a_1=0$. Then the left side of 
(\ref{ZT3.5B}) 
equals $(1+\theta_d)^{-1}$, and is therefore negative. The right side 
of 
(\ref{ZT3.5B}) 
equals
$(1+\theta_1)^{-1}$, and is therefore positive. 
Next assume $a_1\not=0$. Then 
(\ref{ZT3.5B}) 
is immediate from 
(\ref{ZT3.3A}).

\noindent (ii) First assume $\Gamma $ is  bipartite. Then 
$\theta_d=-k$ and  $a_1=0$, so both sides of  (\ref{ZT3.5A}) equal 0.
Next assume $\Gamma $ is not bipartite.
Then (\ref{ZT3.5A}) is 
immediate from (i) above, Lemma 
\ref{ZT3.50},
and Lemma \ref{ZL2.0}.
\qed

%
%

\medskip
\noindent
We now consider when  equality is attained in Theorem \ref{ZT3.5}. To avoid
trivialities, we consider only the nonbipartite case.  
%

\begin{cor}\label{corbf3} \label{NC3.11} \label{ZC3.6} 
Let $\Gamma$ denote a nonbipartite distance-regular graph with diameter 
$d \geq 3$, and  eigenvalues
$\; \theta_0 > \theta_1 > \cdots > \theta_d$. 
 Then the following are equivalent.
\begin{description} \itemsep -3pt
\item{(i)}   \ \ \ Equality holds in  (\ref{ZT3.5B}).
\item{(ii)}   \ \  Equality holds in  (\ref{ZT3.5A}).
\item{(iii)} \ \   $a_1\not= 0$ and every edge of $\Gamma $ 
                    is tight with respect to
                   both $\theta_1$ and $\theta_d$. 
\item{(iv)}   \ \ \   $a_1\not= 0$ and there exists an edge of $\Gamma $
                    which is 
		   tight with respect to both $\theta_1$ and $\theta_d$.
\end{description}
\end{cor}

\proof (i) \iff (ii)
Immediate from Lemma \ref{ZT3.50}. 

\noindent (i),(ii) \implies (iii)
Suppose $a_1=0$. We assume  (\ref{ZT3.5A}) holds with equality, so
$
(\theta_1 + k)(\theta_d + k) = 0,
$
forcing $\theta_d = -k$. Now $\Gamma$ is bipartite by Lemma \ref{drgbip}, 
contradicting the assumption. Hence $a_1 \not= 0$.
Let $xy$ denote an edge of $\Gamma $. Observe the expressions on the left and right
in (\ref{ZT3.3A}) are equal, so they both equal $f(x,y)$.
Now $xy$ is tight with respect to both
$\theta_1$,  $\theta_d$ by Corollary \ref{NT3.10}(i),(ii). 
\\
(iii) \implies (iv) Clear.
\\
(iv) \implies (i) Suppose the edge $xy$ is tight with respect
to both $\theta_1$, $\theta_d$. By Corollary \ref{NT3.10}(i),(ii),
 the scalar
$f(x,y)$ equals both the expression
on the left and the expression on the right in 
 (\ref{ZT3.3A}),  so these expressions are equal.
  \qed

\begin{df} \label{D3.2} \label{AD3.8}                
Let $\Gamma =(X,R)$ denote a distance-regular graph with 
diameter $d\geq 3$.
We say  $\Gamma$ is {\bf tight} whenever 
%
%
$\Gamma$ is not bipartite and the equivalent conditions (i)--(iv) hold in
Corollary \ref{NC3.11}.
\end{df}

\noindent  We wish to emphasize the following fact.
\begin{prop} \label{ZP3.7}                          
Let $\Gamma$ denote a tight distance-regular graph with diameter 
$d \ge 3$. Then $a_i \not= 0$ \\ $(1 \le i \le d-1)$.
\end{prop}

\proof Observe $a_1 \not= 0$ by Corollary \ref{ZC3.6}(iii) and 
 Definition \ref{AD3.8}. Now $a_2,\dots,a_{d-1}$ are nonzero by 
 Lemma~\ref{ZZL2.11}.\qed

\medskip
\noindent We finish this section with some inequalities involving the
eigenvalues of tight graphs.

\begin{lemma} \label{AL3.9}                              
Let $\Gamma = (X,R)$ denote a tight distance-regular graph with
diameter $d \geq 3$ and  eigenvalues 
$\theta_0 > \theta_1 > \cdots > \theta_d$. Then (i)--(iv) hold below.
\begin{description} \itemsep -3pt
\item{(i)}\ \ \ $\theta_d < {\displaystyle{-k \over a_1+1}}$.
\item{(ii)}\ \ Let $\rho$, $\rho_2$ denote the first and second 
	     cosines for $\theta_d$, respectively. Then  
$\;            \rho^2 < \rho_2$.             %
	    \label{AL3.9A}
%
%
\item{(iii)}  Let $\sigma$, $\sigma_2$ denote the first and second 
	     cosines for $\theta_1$, respectively. Then  
         $\;   \sigma^2 > \sigma_2$.             
%
\item{(iv)}\ For each edge $xy $ of $\Gamma$, the scalar $f=f(x,y)$ satisfies
$0<f< b_1$.
	     %
\end{description}
\end{lemma}
 
\proof $\!\!$(i) 
Observe (\ref{ZT3.5A}) holds with equality since $\Gamma $ is tight, and
$a_1 \not= 0$ by Proposition \ref{ZP3.7}, so  
$$
\left(\theta_1 + {k \over a_1+1}\right)\left(\theta_d + {k \over a_1+1}
\right)
< 0.
$$
Since $\theta_1 > \theta_d$, the first factor is positive, and the 
second is negative. The result follows.
\\
(ii) 
By Lemma \ref{TECH}(iv),
\begin{equation}
k^2b_1(\rho^2-\rho_2)
\ = \ 
(k-\theta_d)(k+\theta_d(a_1+1)).     \label{AL3.9C}
\end{equation}
The right side of (\ref{AL3.9C}) is negative in view of (i) above, so
 $\rho^2 < \rho_2$.
\\
(iii) 
By Lemma \ref{TECH}(iv),
\begin{equation}
k^2b_1(\sigma^2-\sigma_2)
\ = \ 
(k-\theta_1)(k+\theta_1(a_1+1)).     \label{AL3.9Z}
\end{equation}
The right side of (\ref{AL3.9Z}) is positive in view of Lemma
\ref{ZL2.0}(i),
 so
 $\sigma^2 > \sigma_2$.
\\
(iv) Observe $f$ equals the expression on the right
in (\ref{ZT3.3A}). This expression is positive and less than $b_1$,
since $\theta_1$ is positive.  \qed

\section{Two characterizations of tight graphs}         

\begin{th}\label{ZT3.11}                           
Let $\Gamma$ denote a nonbipartite distance-regular graph with 
diameter $\;d\ge 3$, and  eigenvalues 
$\; \theta_0 >\theta_1 >\cdots >\theta_d$. Then for all real numbers
$\alpha, \beta$, the following are equivalent.
\begin{description} \itemsep -1pt
\item{(i)} \ \ $\Gamma $ is tight, and $\alpha, \beta $ 
               is a permutation of $\theta_1, \theta_d$.
\item{(ii)} \  $\theta_d \le \alpha, \beta \le  \theta_1$, and 
               \begin{equation}
               \left(\alpha + {k \over a_1+1}\right)
               \left(\beta + {k \over a_1+1} \right)
               = -{ka_1b_1 \over (a_1+1)^2}.         \label{ZT3.11A}
               \end{equation}
\end{description}
\end{th}

\proof
(i) \implies (ii) Immediate since (\ref{ZT3.5A}) holds with equality.
\\
(ii) \implies (i)  Interchanging  $\alpha$ and $\beta$ if necessary, 
                    we may assume  $\alpha \ge \beta$. Since the right 
		    side  of (\ref{ZT3.11A}) is nonpositive, we have
\begin{eqnarray*}
0 \;\leq \; \alpha + {k \over a_1+1} \;&\leq\; \theta_1 
	      +{\displaystyle {k \over a_1+1}}, \\
0 \;\geq \; \beta +  {k \over a_1+1} \;&\geq\; \theta_d 
	      +{\displaystyle {k \over a_1+1}}.
\end{eqnarray*}
By (\ref{ZT3.11A}), the above inequalities, and (\ref{ZT3.5A}), we have 
\begin{eqnarray}
-{ka_1b_1 \over (a_1+1)^2} &=& 
\left(\alpha +{\displaystyle {k \over a_1+1}} \right)
\left(\beta  +{\displaystyle {k \over a_1+1}} \right)
\nonumber
\\
&\geq &
\left(\theta_1 +{\displaystyle {k \over a_1+1}} \right)
\left(\theta_d +{\displaystyle {k \over a_1+1}} \right)  \label{ZT3.11B}
\\
&\geq &  -{\displaystyle{ka_1b_1 \over (a_1+1)^2}}.  \label{ZT3.11C}
\end{eqnarray}
Apparently we have equality in (\ref{ZT3.11B}), (\ref{ZT3.11C}). 
In particular (\ref{ZT3.5A}) holds with equality, so $\Gamma $ is tight.
We mentioned equality holds in (\ref{ZT3.11B}). Neither side is $0$, 
since $a_1\not=0$ by Proposition \ref{ZP3.7}, and it follows 
$\alpha =\theta_1$, $\beta= \theta_d$. \qed
\\


\begin{th}                             \label{ZZT4.3}   
Let $\Gamma=(X,R)$ denote a nonbipartite distance-regular graph with 
diameter $d\ge 3$, and eigenvalues $\theta_0>\theta_1 >\cdots>\theta_d$.
Let $\theta$ and $\theta'$ denote  distinct eigenvalues of $\;\Gamma$, 
with respective cosine sequences $\sigma_0,\sigma_1,\dots,\sigma_d$ and
$\rho_0, \rho_1,\dots,\rho_d$. The following are equivalent.
\begin{description} \itemsep -1pt
\item{(i)}  \ \ \ $\Gamma$ is tight, and $\;\theta$, $\theta'\;$
            is a permutation of $\;\theta_1$, $\theta_d$.
\item{(ii)} \ For $1\leq i \leq d$,
	    \begin{equation}                           \label{ZZT4.3A} 
            {\sigma \sigma_{i-1} - \sigma_i \over
             (1+\sigma) (\sigma_{i-1} - \sigma_i)}
            \;=\; {\rho \rho_{i-1} - \rho_i \over
                  (1+\rho) (\rho_{i-1} - \rho_i)}, 
	    \end{equation}
            and the denominators in (\ref{ZZT4.3A}) are  nonzero. 
\item{(iii)} \begin{equation}                           \label{ZZT4.3B}
            {\sigma^2 - \sigma_2 \over (1+\sigma)(\sigma - \sigma_2)}
            \;=\; {\rho^2 - \rho_2 \over (1+\rho)(\rho - \rho_2)}, 
	    \end{equation}
            and the denominators in (\ref{ZZT4.3B}) are nonzero.
\item{(iv)} \ $\theta$ and $\theta'$ are both nontrivial, and 
            \begin{equation}                  \label{ZZT4.3C} 
            (\sigma_2 \rho_2 - \sigma \rho)(\rho -\sigma)
            \;=\; (\sigma \rho_2 - \sigma_2 \rho)(\sigma \rho -1). 
            \end{equation}
\end{description}
\end{th}

\proof (i) \implies (ii) Recall $a_1\not=0$ by Proposition \ref{ZP3.7}.
Pick adjacent vertices  $x, y \in X$, and let $D^1_1=D_1^1(x,y)$ be
as in Definition \ref{ZD4.1}. By 
 Corollary \ref{ZC3.6}(iii), 
the edge $xy$ is tight with respect to both
$\theta$, $\theta'$; applying (\ref{ZZT4.2A}),
we find both sides of (\ref{ZZT4.3A}) equal 
$a_1^{-1}\vert \Gamma_{i-1}(z)\cap D_1^1 \vert $,
where $z$ denotes any vertex in $D^i_{i-1}(x,y)$.
In particular, the two sides of
 (\ref{ZZT4.3A}) are  equal. 
The denominators 
 in (\ref{ZZT4.3A})  
are nonzero by Lemma \ref{ZL2.6} and
Lemma \ref{drgbip}.
\\
(ii) \implies (iii) Set $i=2$ in (ii). 
\\
(iii) \implies (iv) 
$\theta$ is nontrivial; otherwise $\sigma=\sigma_2=1$, and a denominator
in 
(\ref{ZZT4.3B}) is zero. 
Similarly $\theta'$ is nontrivial. To get (\ref{ZZT4.3C}),  put 
(\ref{ZZT4.3B}) over a common denominator and simplify the result.
\\
(iv) \implies (i) Eliminating $\sigma, \sigma_2, \rho, \rho_2$ in
(\ref{ZZT4.3C}) using $\theta = k\sigma$, $\theta'=k \rho$,  
and Lemma \ref{TECH}(i), 
we routinely find 
 (\ref{ZT3.11A}) holds for $\alpha=\theta$ and $\beta=\theta'$.
Applying Theorem \ref{ZT3.11}, we find $\Gamma $ is tight, and that
 $\theta$, $\theta'$ is a permutation of
$\theta_1$, $\theta_d$.
 \qed

\section{The auxiliary parameter}       

Let $\Gamma $ denote a tight distance-regular graph with diameter 
$d\ge 3$. We are going to show the intersection numbers of $\Gamma$ are
given by certain  rational expressions involving $d$ independent 
parameters. We begin by introducing one of these parameters.

\begin{df} \label{aux}                            
Let $\Gamma$ denote a tight distance-regular graph with diameter 
$d \ge 3$, and  eigenvalues
$\theta_0 > \theta_1 > \cdots >\theta_d$.
Let $\theta$ denote one of  $\theta_1$, $\theta_d$.
By the {\de auxiliary parameter} of $\Gamma$ associated with $\theta$,
we mean the scalar 
\begin{equation}  \label{ZZT4.4A}
\eps \;=\; {k^2-\theta \theta' \over k(\theta - \theta')},
\end{equation}
where $\theta'$ denotes the complement of $\theta $ in 
$\lbrace \theta_1, \theta_d\rbrace$.  We observe
the auxiliary parameter for $\theta_d$ is the opposite of the auxiliary
parameter for $\theta_1$.
\end{df}

\begin{lemma} \label{auxbound}                            
Let $\Gamma$ denote a tight distance-regular graph with diameter 
$d \ge 3$, and  eigenvalues
$\theta_0 > \theta_1 > \cdots >\theta_d$. 
Let $\theta $ denote one of 
$ \theta_1, \theta_d$, and let $\eps $ denote
the auxiliary parameter for $\theta $. Then (i)--(iv) hold below.
\begin{description} \itemsep -1pt
\item{(i)} \ \ \ $\eps > 0 \;$  if $\;\theta=\theta_1$, and $\;\eps < 0 $ if
$\theta=\theta_d$.
\item{(ii)}\ \ $1 < |\eps |$. 
\item{(iii)}\ $|\eps | < k\theta^{-1}_1$.
\item{(iv)}   
\ $|\eps | < -k\theta^{-1}_d$.
\end{description} 
\end{lemma}

\proof
First assume $\theta=\theta_1$.  By (\ref{ZZT4.4A}),  
$$
\eps -1 \;=\; (k+\theta_d)(k-\theta_1)(\theta_1-\theta_d)^{-1}k^{-1} \;>\; 0,
$$
so $\eps > 1$. Recall
$\theta_1>0$ and $\theta_d<0$.
By this and 
 (\ref{ZZT4.4A}),  
$$
k\theta_1^{-1}
-\eps   \;=\; 
\theta_d(k-\theta_1)(k+\theta_1)(\theta_d-\theta_1)^{-1}k^{-1}\theta_1^{-1} \;>\; 0,
$$
so $\eps < k\theta^{-1}_1 $. Similarily 
$$
k\theta^{-1}_d+\eps  \;=\; 
\theta_1(k-\theta_d)(k+\theta_d)(\theta_1-\theta_d)^{-1}k^{-1}\theta_d^{-1} \;<\; 0,
$$
so $\eps < -k\theta_d^{-1} $.  
We now have the result for $\theta = \theta_1$. The result for $\theta=\theta_d$
follows in view of the  last line of Definition \ref{aux}. 
\qed

\begin{th} \label{ZZT4.4}                               
Let $\Gamma$ denote a nonbipartite distance-regular graph with diameter $d \ge 3$,
and   eigenvalues $\theta_0 > \theta_1 > \cdots >\theta_d$.
Let $\theta$ and $\theta'$ denote any eigenvalues of 
$\;\Gamma$, with respective cosine sequences 
$\sigma_0,\sigma_1,\dots,\sigma_d$ and $\rho_0, \rho_1,\dots,\rho_d$. 
Let $\eps$ denote any complex scalar.  Then the  following are equivalent.
\begin{description} \itemsep -1pt
\item{(i)}  \ \ \  $\Gamma$ is tight, $\,\theta, \theta'\,$
             is a permutation of $\,\theta_1,\theta_d$, and
	     $\eps $ is the auxiliary parameter for $\theta$.
\item{(ii)} \ \  $\theta$ and $\theta'$ are both nontrivial, and  
             \begin{equation}
	     \sigma_i\rho_i - \sigma_{i-1}\rho_{i-1} \;=\; 
	     \eps (\sigma_{i-1}\rho_i - \rho_{i-1}\sigma_i)
	        	                         \label{ZZT4.4B}
             \end{equation}
             \ for  $1 \leq i \leq d$.
\item{(iii)} \  $\theta$ and $\theta'$ are both nontrivial, and 
             \begin{equation}
	     \sigma\rho - 1 \;=\; \eps (\rho - \sigma),
	     \ \ \ \ \ \ \ \ \ \ \ \ \ \ \ \ \
	     \sigma_2\rho_2 - \sigma\rho \;=\; 
	     \eps (\sigma\rho_2 - \rho\sigma_2).
	     %
	        	                        \label{ZZT4.4C}
             \end{equation}
\end{description}
\end{th}

\proof 
(i) \implies (ii) 
It is clear $\theta$, $\theta'$ are both nontrivial.
To see (\ref{ZZT4.4B}), observe $\theta, \theta'$ are distinct, so
the equivalent statements (i)--(iv) 
in Theorem \ref{ZZT4.3} hold. Putting (\ref{ZZT4.3A}) over a common
denominator and simplifying using $\eps = (1-\sigma \rho)(\sigma-\rho)^{-1}$,
we get (\ref{ZZT4.4B}).
\\
(ii) \implies (iii) 
%
Set $i=1$ and 
$i=2$ in  
 (\ref{ZZT4.4B}). 
\\
(iii) \implies (i) 
We first show $\theta \not= \theta'$. Suppose $\theta = \theta'$. Then
$\sigma =\rho$, so the left equation of 
(\ref{ZZT4.4C}) becomes $\sigma^2 = 1$, forcing $\sigma = 1$ or 
$\sigma = -1$. But $\sigma\not=1$ since  $\theta$ is 
nontrivial, and $\sigma \not= -1$ since $\Gamma$ is not bipartite. 
We conclude $\theta \not= \theta'$. 
Now $\sigma \not= \rho$; 
solving the left equation in (\ref{ZZT4.4C})
for $\eps$, and 
eliminating $\eps$ in the right equation of
(\ref{ZZT4.4C}) using the result, we obtain (\ref{ZZT4.3C}). Now Theorem 
\ref{ZZT4.3}(iv) holds. 
Applying Theorem \ref{ZZT4.3}, we find $\Gamma$ is tight, and that 
$\theta$, $\theta'$ is a permutation of $\theta_1$, $\theta_d$.
Solving the left equation in (\ref{ZZT4.4C})
for $\eps$, and simplifying the result,
we obtain (\ref{ZZT4.4A}). It follows $\eps $ is the auxiliary 
parameter for $\theta$.  \qed

\section{Feasibility}                                     

Let $\Gamma$ denote a tight distance-regular graph with diameter 
$d \ge 3$, and  eigenvalues
$\theta_0 > \theta_1 > \cdots >\theta_d$. Let
 $\theta,\theta'$ denote a permutation of $\theta_1, \theta_d$,  
 with respective cosine sequences 
$\sigma_0,\sigma_1,\dots,\sigma_d$ and $\rho_0, \rho_1,\dots,\rho_d$. 
Let $\eps$ denote the auxiliary parameter for $\theta $. Pick any
integer $i$ $(1 \leq i \leq d)$, and 
observe   (\ref{ZZT4.4B}) holds. Rearranging terms in that equation,
we find
\begin{equation}
\rho_i     (\sigma_i - \eps \sigma_{i-1}) \;=\; 
\rho_{i-1} (\sigma_{i-1} - \eps \sigma_i).
\label{takealook}
\end{equation}
We would like to solve (\ref{takealook}) for $\rho_i$, but conceivably 
 $\sigma_i - \eps \sigma_{i-1} = 0$. In this section we investigate this 
possibility.

\begin{lemma} \label{ZZL4.5}                                 
Let $\Gamma$ denote a tight distance-regular graph with diameter 
$d \ge 3$, and  eigenvalues
$\theta_0 > \theta_1 > \cdots >\theta_d$. Let
 $\theta,\theta'$ denote a permutation of $\theta_1, \theta_d$,  
with respective cosine sequences  
$\sigma_0,\sigma_1,\dots,\sigma_d$ and $\rho_0, \rho_1,\dots,\rho_d$.
Let $\eps$ denote the auxiliary parameter for $\theta $.
Then for each integer $\,i$ ($1\le i \le d-1$), the following are equivalent:
(i) $\sigma_{i-1} = \eps \sigma_i$;
(ii) $\sigma_{i+1} = \eps \sigma_i$;
(iii) $\sigma_{i-1} = \sigma_{i+1}$;
(iv) $\rho_i = 0$.
Moreover, suppose (i)--(iv) hold. Then $\theta = \theta_d$ 
and $\theta' = \theta_1$. 
\end{lemma}
%

\proof Observe Theorem \ref{ZZT4.4}(i) holds, so 
(\ref{ZZT4.4B}) holds.

\noindent (i) \implies (iv) 
Replacing $\sigma_{i-1}$ by $\eps\sigma_i$ in
(\ref{ZZT4.4B}),
 we find
$
\sigma_i\rho_i (1-\eps^2) \;=\; 0.
$
Observe $\eps^2\not=1$ by Lemma \ref{auxbound}(ii). Suppose for the 
moment that 
$\sigma_i = 0$. We assume  $\sigma_{i-1} = \eps\sigma_i$, so 
$\sigma_{i-1}=0$. Now $\sigma_{i-1}=\sigma_i$, contradicting
 Lemma \ref{ZL2.6}. Hence $\sigma_i\not=0$, 
so $\rho_i=0$.
\\
(iv) \implies (i) 
Setting $\rho_i=0$ in (\ref{ZZT4.4B}), we find
$
\rho_{i-1} (\sigma_{i-1} - \eps\sigma_i) \;=\; 0.
$
Observe $\rho_{i-1} \not= 0$, otherwise $\rho_{i-1} = \rho_i$,
contradicting Lemma \ref{ZL2.6}. We conclude $\sigma_{i-1}= \eps
\sigma_i$, as desired.
\\
(ii) \iff (iv) Similar to the proof of (i) \iff (iv). 
\\
(i),(ii) \implies (iii) Clear.
\\
(iii) \implies (i)  We cannot have $\theta=\theta_1$ by Lemma \ref{ZL2.6}(i), 
so $\theta =\theta_d$, $\theta'=\theta_1$. In particular 
$\;\rho_{i-1}\not=\rho_{i+1}$.
Adding (\ref{ZZT4.4B}) at $i$ and $i+1$, we obtain
$$
\sigma_{i+1}\rho_{i+1} - \sigma_{i-1}\rho_{i-1}
\;=\; \eps (
\sigma_i\rho_{i+1} - \sigma_{i+1}\rho_i + 
\sigma_{i-1}\rho_i - \sigma_i\rho_{i-1}).
$$
Replacing $\sigma_{i+1}$ by $\sigma_{i-1}$ in the above line,
and simplifying, we obtain
$$
(\sigma_{i-1} - \eps\sigma_i) (\rho_{i+1} - \rho_{i-1}) \;=\; 0.
$$
It follows 
$\sigma_{i-1} = \eps\sigma_i$, as desired.
\\
Now suppose (i)--(iv). Then we saw in the proof of (iii) \implies (i) that
$\theta = \theta_d$, $\theta' = \theta_1$.  \qed

\begin{df} \label{ZZD4.7}                            
Let $\Gamma = (X, R)$ denote a tight distance-regular graph with 
diameter $d \ge 3$ and
eigenvalues $\theta_0 > \theta_1 >\cdots >\theta_d$.
Let $\sigma_0,\sigma_1,\dots,\sigma_d$ denote any cosine sequence for
$\Gamma$ and let $\theta$ denote the corresponding eigenvalue. The 
sequence $\sigma_0,\sigma_1,\ldots,\sigma_d$ (or $\theta$)
is said to be {\de feasible} whenever (i) and (ii) hold below.
\begin{description} \itemsep -1pt
\item{(i)}  \ \ $\theta $ is one of $\,\theta_1, \theta_d$. 
\item{(ii)}  \    
                $\sigma_{i-1} \not= \sigma_{i+1}$ for   $\,1 \le i \le d-1$.
\end{description} 
We observe by Lemma  \ref{L5.1}(i)
that $\theta_1$ is feasible.
\end{df}
We conclude this section with an extension of Theorem \ref{ZZT4.4}.

\begin{th} \label{ZZL4.6}                               
Let $\Gamma$ denote a nonbipartite
distance-regular graph with diameter $d \ge 3$,
 and  eigenvalues $\theta_0 > \theta_1 > \cdots >\theta_d$.
Let $\theta$ and $\theta'$ denote any  eigenvalues of $\;\Gamma$,
with respective cosine sequences 
$\sigma_0,\sigma_1,\dots,\sigma_d$ and $\rho_0, \rho_1,\dots,\rho_d$. 
Let $\eps$ denote any complex scalar.
Then the following are equivalent.
\begin{description} \itemsep -1pt
\item{(i)}  \ \ $\Gamma $ is tight, $\theta $ is feasible,  
		 $\eps$  is the auxiliary
		 parameter for $\theta$,
and $\theta'$ is the complement of $\theta $ in $\lbrace \theta_1, \theta_d \rbrace $.  
\item{(ii)} \   $\theta'$ is not trivial, 
                \begin{equation}
                \rho_i \;=\; 
		\prod_{j=1}^i {\sigma_{j-1} - \eps \sigma_j \over
                               \sigma_j - \eps \sigma_{j-1}}
                               \qquad \qquad  (0 \le i \le d),
	        	                        \label{ZZL4.6A}
                \end{equation}
                and denominators in (\ref{ZZL4.6A}) are all nonzero.
\end{description}
\end{th}

\proof
(i) \implies (ii) Clearly $\theta'$ is nontrivial.
To see (\ref{ZZL4.6A}), observe 
Theorem \ref{ZZT4.4}(i) holds, so 
(\ref{ZZT4.4B}) holds.
Rearranging  terms in (\ref{ZZT4.4B}), we obtain
\begin{equation}
\rho_i (\sigma_i - \eps \sigma_{i-1}) \;=\; 
\rho_{i-1} (\sigma_{i-1} - \eps \sigma_i) \qquad \qquad (1 \le i \le d).
\label{ZZL4.6B}
\end{equation}
Observe $\sigma_i \not= \eps \sigma_{i-1}$ for $2 \le i \le d$ by 
Lemma  \ref{ZZL4.5}(ii), and $\sigma \not = \eps $ by Lemma \ref{auxbound}(ii),
so the coefficient of $\rho_i$ in (\ref{ZZL4.6B}) is never zero. Solving 
that equation for $\rho_i$ and applying induction, we  routinely obtain 
(\ref{ZZL4.6A}).
\\
(ii) \implies (i) We show Theorem \ref{ZZT4.4}(iii) holds.
Observe $\theta$ is nontrivial; otherwise $\sigma =1$, forcing 
$\rho =1$ by (\ref{ZZL4.6A}), and contradicting our assumption that 
$\theta'$ is nontrivial. One readily verifies (\ref{ZZT4.4C}) by
eliminating $\rho, \rho_2$ using (\ref{ZZL4.6A}).
We now 
have Theorem \ref{ZZT4.4}(iii). Applying that theorem,
we find
 $\Gamma $ is tight, $\theta, \theta'\,$ is a permutation
 of $\theta_1, \theta_d$, and that  $\eps$  is the auxiliary
 parameter for $\theta$.
It remains to show $\theta $ is feasible. Suppose not.
Then there exists an integer $i $ $(1 \leq i \leq d-1)$ such
that $\sigma_{i-1}= \sigma_{i+1}$. Applying Lemma \ref{ZZL4.5}, we
find $\sigma_{i+1}=\eps \sigma_i$. But 
$\sigma_{i+1}-\eps \sigma_i$ is a  factor in the denominator of
(\ref{ZZL4.6A}) (with $i$ replaced by $i+1$), and hence is not 0.
We now have a contradiction, so $\theta $ is feasible.   \qed

\section{A parametrization}                        

In this section, we obtain the intersection numbers
of a tight graph as rational functions of a feasible cosine
sequence and the associated auxiliary parameter.
We begin
with a  result about arbitrary distance-regular graphs.

\begin{lemma} \label{ZT2.9}          
Let $\Gamma$ denote  a distance-regular graph with diameter $d\ge 3$,
 and eigenvalues  $\theta_0 > \theta_1 >\cdots >\theta_d$.
Let $\theta$, $\theta'$ denote a permutation of $\theta_1$, $\theta_d$,
with respective cosine sequences $\;\sigma_0,\sigma_1,\ldots,\sigma_d$
and $\;\rho_0,\rho_1,\ldots,\rho_d$. Then
\begin{eqnarray}
k \;&=& \; {\displaystyle{
       {(\sigma-\sigma_2)(1-\rho)-(\rho-\rho_2)(1-\sigma) }
       \over 
       {(\rho-\rho_2)(1-\sigma)\sigma-(\sigma-\sigma_2)(1-\rho)\rho}}}, 
				        \label{ZT2.9A}\\
b_i\;&=& \; k{\displaystyle{
	  {(\sigma_{i-1}-\sigma_i)(1-\rho)\rho_i
          - (\rho_{i-1}-\rho_i)(1-\sigma)\sigma_i }
          \over 
          {(\rho_i-\rho_{i+1})(\sigma_{i-1}-\sigma_i)
          -(\sigma_i-\sigma_{i+1})(\rho_{i-1}-\rho_i)}}} 
	  \quad (1\leq i \leq d-1), \qquad \qquad       
				        \label{ZT2.9B}\\
c_i\;&=&\;k{\displaystyle{
	  {(\sigma_i-\sigma_{i+1})(1-\rho)\rho_i
	  -(\rho_i-\rho_{i+1})(1-\sigma)\sigma_i    }
	  \over 
	  {(\rho_i-\rho_{i+1})(\sigma_{i-1}-\sigma_i)
	  -(\sigma_i-\sigma_{i+1})(\rho_{i-1}-\rho_i)}}} 
	  \quad (1\leq i \leq d-1),                     
					\label{ZT2.9C}\\
c_d \;&=&\; k\sigma_d {\sigma-1 \over \sigma_{d-1}-\sigma_d}
     \;=\; k\rho_d   {\rho-1   \over \rho_{d-1}-\rho_d},
						 \label{ZT2.9D}
\end{eqnarray}
and the denominators in (\ref{ZT2.9A})--(\ref{ZT2.9D}) are never zero. 
\end{lemma}

\proof Line (\ref{ZT2.9D}) is immediate from Lemma \ref{TECH}(v),
and the denominators in that line are nonzero by Lemma \ref{ZL2.6}.
To obtain (\ref{ZT2.9B}), 
 (\ref{ZT2.9C}), pick any integer $i$ $(1 \leq i \leq d-1)$, and recall  
by Lemma \ref{ZL2.4}(iii) that
\begin{eqnarray}
  c_i(\sigma_{i-1}-\sigma_i)-b_i(\sigma_i-\sigma_{i+1})
  &=& k(\sigma-1)\sigma_i,             \label{ZT2.9E}
\\
  c_i(\rho_{i-1}-\rho_i)-b_i(\rho_i-\rho_{i+1})
  &=& k(\rho-1)\rho_i.                 \label{ZT2.9F}
\end{eqnarray}
To solve this linear system for $\,c_i$ and $\,b_i$,
 consider the determinant
$$ 
D_i := \hbox{det} \left( \begin{array}{cc} 
    \sigma_{i-1}-\sigma_i & \sigma_i-\sigma_{i+1}\cr
    \rho_{i-1}-\rho_i     & \rho_i-\rho_{i+1}    \cr
    \end{array}\right).                             
$$
Using  Lemma \ref{ZL2.6}, we routinely find
$\;D_i\not=0$.
Now (\ref{ZT2.9E}), (\ref{ZT2.9F}) has the unique solution 
(\ref{ZT2.9B}), (\ref{ZT2.9C}) by elementary linear algebra. 
The denominators in 
(\ref{ZT2.9B}), (\ref{ZT2.9C}) 
both  equal $D_i$;
in particular they 
are not zero. To get (\ref{ZT2.9A}), set $\;i=1\; $
and $\;c_1=1\;$ in (\ref{ZT2.9C}),
and solve for~$\;k$. 
 \qed

\begin{th} \label{ZZT4.8}                             
Let $\Gamma$ denote a nonbipartite 
distance-regular graph with diameter $d \ge 3$, and let  
 $\;\sigma_0, \sigma_1, \ldots \sigma_d,\eps, h\;$ denote
complex scalars. Then the following are equivalent.
\begin{description}
\item{(i)} \ \  $\Gamma$ is tight, $\sigma_0, \sigma_1, \ldots \sigma_d$
                is a feasible cosine sequence for $\Gamma$, $\eps$ is 
		the associated auxiliary parameter from (\ref{ZZT4.4A}),
		and 
                \begin{equation}
                h \;=\; {{(1-\sigma)(1-\sigma_2)} \over 
                        {(\sigma^2-\sigma_2)(1 -\eps\sigma)}}.  
			\label{ZZT4.8A}
                \end{equation}
\item{(ii)} \  $\sigma_0=1$, $\sigma_{d-1} = \sigma \sigma_d$,
               $\eps \not= -1$, 
\begin{eqnarray}
k \;&=&\;  {\displaystyle
	   h {\sigma -\eps \over \sigma - 1}},  \label{ZZT4.8B}
\\[1ex]
b_i \;&=&\; {\displaystyle 
	    h{(\sigma_{i-1}-\sigma \sigma_i)(\sigma_{i+1}-\eps \sigma_i)
              \over (\sigma_{i-1}-\sigma_{i+1})(\sigma_{i+1}-\sigma_i)}}
	           \qquad (1 \leq i \leq d-1),\label{ZZT4.8C}
\\[1ex]
c_i \;&=&\; {\displaystyle 
	    h{(\sigma_{i+1}-\sigma \sigma_i)(\sigma_{i-1}-\eps \sigma_i)
              \over (\sigma_{i+1}-\sigma_{i-1})(\sigma_{i-1}-\sigma_i)}}
	           \qquad (1 \leq i \leq d-1),\label{ZZT4.8D}
\\[1ex]
c_d \;&=&\;  {\displaystyle 
	   h{\sigma -\eps \over \sigma - 1}},\label{ZZT4.8E}
\end{eqnarray}
and denominators in (\ref{ZZT4.8B})--(\ref{ZZT4.8E}) are all nonzero.
\end{description}
\end{th}

%
%

\proof Let $\theta_0> \theta_1> \cdots > \theta_d$ denote the eigenvalues
of $\Gamma $.

\noindent (i) \implies (ii)
%
Observe $\sigma_0=1$ by Lemma \ref{ZL2.4}(ii), and 
$\eps \not= -1$ 
 by Lemma  \ref{auxbound}(ii). 
Let $\theta$ denote the eigenvalue associated with 
$\sigma_0,\sigma_1,\dots,\sigma_d$, and observe by 
Definition 
 \ref{ZZD4.7} that   
$\theta$ is one of 
$\theta_1,\theta_d$. Let $\theta'$ denote the complement of 
$\theta$
in $\{ \theta_1,\theta_d\}$,  and let 
$\rho_0, \rho_1,\dots,\rho_d$ denote the cosine sequence for $\theta'$.
%
%
Observe  Theorem  \ref{ZZL4.6}(i) holds. Applying that 
theorem, we obtain (\ref{ZZL4.6A}).
%
%
Eliminating 
$\rho_0, \rho_1,\dots,\rho_d$ in (\ref{ZT2.9A})--(\ref{ZT2.9D}) using
(\ref{ZZL4.6A}), we routinely obtain
(\ref{ZZT4.8B})--(\ref{ZZT4.8E}),
and that $\sigma_{d-1} = \sigma\sigma_d$.
\\
(ii) \implies (i)
One readily checks
$$
c_i(\sigma_{i-1} - \sigma_i) - b_i(\sigma_i-\sigma_{i+1}) 
= k(\sigma-1) \sigma_i \qquad \qquad (1 \le i \le d),
$$
where  $\sigma_{d+1}$ is an indeterminant.  Applying   
Lemma \ref{ZL2.4}(i),(iii), we find  $\sigma_0,\sigma_1,\dots,\sigma_d$ is 
a cosine sequence for $\Gamma$, with associated  eigenvalue $\theta := k\sigma$.
By (\ref{ZZT4.8B}), (\ref{ZZT4.8C}), and since $k, b_1,\ldots,b_{d-1}\;$ are
nonzero,
$$
\sigma_j \not= \eps \sigma_{j-1} \qquad \qquad (1 \le j \le d).
$$
Set
\begin{equation}
     \rho_i := \prod_{j=1}^i {\sigma_{j-1} - \eps \sigma_j     \over
                              \sigma_j     - \eps \sigma_{j-1}}
                                    \qquad \qquad  (0 \le i \le d).
	        	                        \label{ZZT4.8G}
\end{equation}
One readily checks $\rho_0 = 1$, and  that
$$
c_i (\rho_{i-1} -\rho_i) -b_i (\rho_i -\rho_{i+1}) \;=\; k (\rho-1)\rho_i
                                            \qquad \qquad  (1 \le i \le d),
$$
where  $\rho_{d+1}$ is an indeterminant. Applying   
Lemma \ref{ZL2.4}(i),(iii), we find  $\;\rho_0,\rho_1,\ldots,\rho_d\;$ is a 
cosine sequence for $\Gamma$, with associated eigenvalue $\theta' := k\rho$.
We claim $\theta'$ is not trivial. Suppose $\theta'$ is trivial. 
Then $\rho = 1$. Setting $i=1$ and $\rho =1$ in (\ref{ZZT4.8G}) we
find $\;\sigma -\eps = 1-\eps \sigma\;$, forcing $\;(1-\sigma)(1+\eps)=0$.
Observe $\sigma \not= 1$ since  the
denominator in (\ref{ZZT4.8E}) is not zero, and we assume $\eps \not= -1$, 
so we have a contradiction. We have now shown $\theta'$ is 
nontrivial, so  
%
%
%
%
%
%
Theorem \ref{ZZL4.6}(ii) holds.
Applying that theorem, 
%
we find $\Gamma$ is tight, $\theta$ is feasible, and 
that $\eps$ is the auxiliary
parameter of $\theta$.
%
%
To see (\ref{ZZT4.8A}), set $i=1$ and $c_1=1$ in (\ref{ZZT4.8D}),
and solve for
$h$.  \qed

\begin{prop}\label{ZZL4.8A}                              
With the notation of Theorem \ref{ZZT4.8},
suppose (i), (ii) hold, and 
let $\theta_0>\theta_1>\cdots
> \theta_d$ denote the  eigenvalues of $\Gamma $.
If $\eps > 0$, then 
\begin{eqnarray}
\theta_1 \;= \: {\displaystyle 
		   {\sigma (\sigma-\eps)(1-\sigma_2) \over
		  (1-\eps\sigma)(\sigma_2 - \sigma^2)}}, 
\qquad \qquad 
\theta_d \;= \: {\displaystyle 
		  {1-\sigma_2 \over \sigma_2 - \sigma^2}}.
                                               \label{ZZL4.8AB} 
\end{eqnarray}
If $\eps < 0$, then 
\begin{eqnarray}
\theta_1 \;=\:  {\displaystyle{1-\sigma_2 \over \sigma_2 - \sigma^2}},
\qquad \qquad 
\theta_d \;=\:  {\displaystyle {\sigma (\sigma-\eps)(1-\sigma_2) \over
		  (1-\eps\sigma)(\sigma_2 - \sigma^2)}}. 
                                               \label{ZZL4.8AD} 
\end{eqnarray}
We remark that the denominators in (\ref{ZZL4.8AB}), (\ref{ZZL4.8AD})
are nonzero.
\end{prop}

\proof Let $\theta $ denote the eigenvalue of $\Gamma $ associated with
$\sigma_0, \sigma_1, \ldots, \sigma_d$. By
Lemma \ref{ZL2.4}(iii) and (\ref{ZZT4.8B}),
we obtain  
\begin{eqnarray}
\theta &=& k \sigma  \nonumber
\\
&=& {\sigma (\sigma-\eps)(1-\sigma_2) \over
		    (1-\eps\sigma)(\sigma_2 - \sigma^2)}. 
\end{eqnarray}
Observe  
$\theta \in \lbrace\theta_1, \theta_d\rbrace $ since
$\;\sigma_0, \sigma_1, \ldots, \sigma_d\;$
is feasible.   
Let $\theta'$ denote the complement of $\theta $ in  
$\lbrace \theta_1, \theta_d\rbrace $, and
let 
$\rho $ denote the first cosine associated with $\theta'$. 
Observe
condition (i) holds in Theorem  \ref{ZZL4.6}, so 
(\ref{ZZL4.6A}) holds. Setting  
 $i=1$ in that equation, 
 we find  
\begin{equation}
\rho= {{1-\eps  \sigma}\over {\sigma -\eps}}. 
\label{rhofromeps}
\end{equation}
 By Lemma \ref{ZL2.4}(iii), (\ref{ZZT4.8B}), and (\ref{rhofromeps}),
we obtain  
\begin{eqnarray}
\theta' &=& k \rho  \nonumber
\\
           &=& {1-\sigma_2 \over
		   \sigma_2 - \sigma^2}. 
\end{eqnarray}
To finish the proof, we observe by  
Lemma \ref{auxbound}(i) that  
 $\theta=\theta_1$,
$\theta'=\theta_d$ 
if $\eps > 0$, and  
$\theta=\theta_d$, 
$\theta'=\theta_1$
if $\eps < 0$.
\qed

\begin{th}\label{ZZT4.9}                              
Let $\Gamma$ denote a tight distance-regular graph with diameter
$\;d\geq 3$, and eigenvalues  $\theta_0 > \theta_1 >\cdots >\theta_d$.
Then (i) and (ii) hold below.
%
\begin{description}\itemsep 30pt \parskip -30pt 
\item{(i)}  \ \ $a_d = 0$.
\item{(ii)} \ Let $\sigma_0,\sigma_1,\dots,\sigma_d$ denote the cosine 
              sequence for $\theta_1$ or $\theta_d$, and let $ \eps$ 
	      denote the associated auxiliary parameter from 
	      (\ref{ZZT4.4A}).
\end{description}
Then 
\begin{equation}
a_i \;=\; g\; {(\sigma_{i+1}-\sigma \sigma_i)(\sigma_{i-1}-\sigma \sigma_i)
     \over (\sigma_{i+1}-\sigma_i)(\sigma_{i-1}-\sigma_i)}
	          \qquad \qquad (1 \leq i \leq d-1),\label{ZZT4.9A}
\end{equation}
where
\begin{equation}
g \;=\; {(\eps-1)(1-\sigma_2) \over (\sigma^2-\sigma_2)(1-\eps\sigma)}.
						    \label{ZZT4.9B}
\end{equation}
\end{th}
 
\proof (i) Comparing   (\ref{ZZT4.8B}), 
 (\ref{ZZT4.8E}), we see   
$k=c_d$,  and it follows $a_d=0$.
\\
\noindent (ii) First assume
$\sigma_0,\sigma_1, \ldots, \sigma_d$ is the cosine sequence for
$\theta_1$, and  recall this sequence is feasible.
 Let $h$ be as in 
(\ref{ZZT4.8A}). 
Then  Theorem \ref{ZZT4.8}(i) holds, so 
 Theorem \ref{ZZT4.8}(ii) holds.
Evaluating the right side of  
$
a_i \;=\; k - b_i - c_i
$
using (\ref{ZZT4.8B})--(\ref{ZZT4.8D}), and simplifying the result using
(\ref{ZZT4.8A}), we obtain (\ref{ZZT4.9A}), (\ref{ZZT4.9B}).
To finish the proof, 
let 
 $\rho_0,\rho_1,\ldots,\rho_d$ denote the cosine sequence for 
$\theta_d$, and recall by Definition \ref{aux} 
that the associated auxiliary parameter is $\,\eps'=-\eps$.
We show
\begin{eqnarray}
                 a_i&=&
 {{(\eps'-1)(1-\rho_2) }\over 
           {(\rho^2-\rho_2)(1-\eps'\rho)}}\,{{
           (\rho_{i+1}-\rho \rho_i)(\rho_{i-1}-\rho \rho_i)}
           \over { 
                 (\rho_{i+1}-\rho_i)(\rho_{i-1}-\rho_i)}}.\label{bigline}
\end{eqnarray}
By Theorem \ref{ZZT4.3}(ii) (with $i$ replaced by $i+1$),
\begin{equation}
{1 \over 1 + \sigma} 
{\sigma_{i+1} - \sigma \sigma_i \over \sigma_{i+1} - \sigma_i}
\;=\; {1 \over 1 + \rho} 
{\rho_{i+1} - \rho \rho_i \over \rho_{i+1} - \rho_i}.
						    \label{ZZT4.9E}
\end{equation}
Subtracting 1 from both sides of Theorem \ref{ZZT4.3}(ii), 
and simplifying, we obtain
\begin{equation}
{1 \over 1 + \sigma} 
{\sigma_{i-1} - \sigma \sigma_i \over \sigma_{i-1} - \sigma_i}
\;=\; {1 \over 1 + \rho} 
{\rho_{i-1} - \rho \rho_i \over \rho_{i-1} - \rho_i}.
						    \label{ZZT4.9F}
\end{equation}
By (\ref{ZZT4.4C}),
\begin{equation}
{(\eps-1)(1-\sigma_2)(1+\sigma)^2 \over 
  (\sigma^2 - \sigma_2)(1- \eps\sigma) } 
\;=\; 
{(\eps'-1)(1-\rho_2)(1+\rho)^2 \over 
  (\rho^2 - \rho_2)(1- \eps'\rho) }.
						    \label{ZZT4.9G}
\end{equation}
Multiplying together 
(\ref{ZZT4.9E})--(\ref{ZZT4.9G}) and simplifying,
we obtain (\ref{bigline}),
as desired.
\qed
\medskip

\noindent We end this section with some inequalities.

\begin{lemma}\label{ZZT4.11}                              
Let $\Gamma$ denote a tight distance-regular graph with diameter
$\;d\geq 3$, and eigenvalues  $\theta_0 > \theta_1 >\cdots >\theta_d$.
Let  $\theta $ denote one of $\theta_1, \theta_d$, and  let 
$\sigma_0,\sigma_1,\ldots, \sigma_d$ denote the  cosine
sequence for $\theta$.

\noindent Suppose $\theta = \theta_1$. Then
%
\begin{description}\itemsep 30pt \parskip -30pt 
\item{(i)}  \ \ $\sigma_{i-1} >  \sigma \sigma_i \qquad \qquad
(1 \leq i \leq d-1), $
\item{(ii)}  \ \ $\sigma \sigma_{i-1} >   \sigma_i \qquad \qquad 
(2 \leq i \leq d)$.
\end{description}
\noindent Suppose $\theta = \theta_d$. Then
%
\begin{description}\itemsep 30pt \parskip -30pt 
\item{(iii)}  \ \ $(-1)^{i}(\sigma\sigma_i -   \sigma_{i-1})>0 \qquad \qquad
(1 \leq i \leq d-1), $
\item{(iv)}  \ \ 
 $(-1)^{i}(\sigma_i -   \sigma\sigma_{i-1})>0 \qquad \qquad
(2 \leq i \leq d)$.
\end{description}
\end{lemma}

\proof (i) We first show $\sigma_{i-1}-\sigma \sigma_i$ is nonnegative.
Recall $a_1\not=0$ by Proposition \ref{ZP3.7}, so Theorem \ref{ZZT4.2}
applies. Let $x,y $ denote adjacent vertices in $X$, and recall by
Corollary \ref{ZC3.6} that the edge $xy$ is tight with respect to
$\theta$. Now Theorem \ref{ZZT4.2}(i) holds, so (\ref{ZZT4.2B}) holds.
Observe 
the left side of  (\ref{ZZT4.2B}) is nonnegative, so the right side
is nonnegative. In that expression on the right, the factors  
$1+\sigma $ and $\sigma_{i-1}-\sigma_i$ are positive, so the remaining factor
$\sigma_{i-1}-\sigma \sigma_i$  is nonnegative, as desired.
To finish the proof, observe  $\sigma_{i-1}-\sigma \sigma_i$ is a factor 
on the right in (\ref{ZZT4.9A}), so it is not zero in view of Proposition
\ref{ZP3.7}.


\noindent (ii)--(iv) Similar to the proof of (i) above.
\qed
\section{The 1-homogeneous property}       

In this section, we show the concept of  tight is closely related to 
the concept of  1-homogeneous that appears in the work of K. Nomura 
\cite{Nom1}, \cite{Nom2}, \cite{Nom3}.

\begin{th} \label{NT5.5} \label{NNL5.8}                 
Let $\Gamma = (X,R)$ denote a tight distance-regular graph with diameter
$d \ge 3$, and eigenvalues $\theta_0 > \theta_1 >\cdots >\theta_d$.
Let $\sigma_0,\sigma_1,\dots,\sigma_d$ denote the cosine sequence 
associated with $\theta_1$ or $\theta_d$.
Fix adjacent vertices $x,y \in X$. 
Then with the notation of Definition \ref{ZD4.1} we have the following:
For all integers $i$ ($1\le i \le d-1$), and for all vertices 
$z \in D_i^i$,
\begin{eqnarray}
|\Gamma_{i-1}(z) \cap D_1^1| \;&=&\; c_i 
       {\displaystyle
       {(\sigma^2 -\sigma_2)(\sigma_i-\sigma_{i+1}) \over
	(\sigma-\sigma_2)(\sigma\sigma_i-\sigma_{i+1})}},
				               \label{NT5.5A}
\\[1ex]
|\Gamma_{i+1}(z) \cap D_1^1| \;&=&\; b_i 
       {\displaystyle
       {(\sigma^2 -\sigma_2)(\sigma_{i-1}-\sigma_i) \over
	(\sigma-\sigma_2)(\sigma_{i-1} - \sigma\sigma_i)}}. 
                                               \label{NT5.5B}
\end{eqnarray}
\end{th}

\proof First assume 
 $\sigma_0,\sigma_1, \ldots, \sigma_d$ is the  cosine sequence for
$\theta_1$, and let $\rho_0,\rho_1,\ldots,\rho_d$ denote the cosine
sequence for $\theta_d$.
%
%
The edge $xy$ is tight with respect to both $\,\theta_1$, $\,\theta_d$,
 so by
Theorem \ref{UZT4.3}(ii), 
\begin{eqnarray}
|\Gamma_{i+1}(z) \cap D_1^1| \;&=&\; |\Gamma_{i-1}(z) \cap D_1^1|\;
{\sigma_{i-1} - \sigma_i \over \sigma_i - \sigma_{i+1}}  
\;+\; a_1\; {1-\sigma \over 1+\sigma} \;
{\sigma_i \over \sigma_i - \sigma_{i+1}},      \label{NT5.5G}
\\[1ex]
|\Gamma_{i+1}(z) \cap D_1^1| \;&=&\; |\Gamma_{i-1}(z) \cap D_1^1|\;
{\rho_{i-1} - \rho_i \over \rho_i - \rho_{i+1}}  
\;+\; a_1\; {1-\rho \over 1+\rho} \;
{\rho_i \over \rho_i - \rho_{i+1}}.                \label{NT5.5H}
\end{eqnarray}
Eliminating $\;\rho_0,\,\rho_1,\,\ldots,\,\rho_d$ in (\ref{NT5.5H})
using (\ref{ZZL4.6A}), we obtain
\begin{eqnarray}
|\Gamma_{i+1}(z) \cap D_1^1| \;&=&\; |\Gamma_{i-1}(z) \cap D_1^1| 
\; {\displaystyle
{\sigma_{i-1} - \sigma_i \over \sigma_i - \sigma_{i+1}}
\; {\sigma_{i+1}-\eps\sigma_i \over\sigma_{i-1}-\eps\sigma_i}}\nonumber 
\\
   &&\;+\;\; a_1 \; {\displaystyle
   {(1-\sigma) (\sigma_{i+1} - \eps\sigma_i) \over
          (1+\sigma)(1-\eps)(\sigma_i-\sigma_{i+1})}}, \label{NT5.5I}
\end{eqnarray}
where  $\eps$ denotes the auxiliary parameter associated with $\theta_1$. 
Solving  (\ref{NT5.5G}), (\ref{NT5.5I})  for 
$|\Gamma_{i+1}(z) \cap D_1^1|$ and 
$|\Gamma_{i-1}(z) \cap D_1^1|$, 
and evaluating the result using 
(\ref{ZZT4.8A}),  
(\ref{ZZT4.8C}), (\ref{ZZT4.8D}), (\ref{ZZT4.9A}), we
get 
(\ref{NT5.5A}), (\ref{NT5.5B}), as desired.  
To finish the proof observe by Theorem \ref{ZZT4.3}(ii),(iii) that
\begin{eqnarray}
       {\displaystyle
       {(\sigma^2 -\sigma_2)(\sigma_i-\sigma_{i+1}) \over
	(\sigma-\sigma_2)(\sigma\sigma_i-\sigma_{i+1})}} 
&=& 		            
       {\displaystyle
       {(\rho^2 -\rho_2)(\rho_i-\rho_{i+1}) \over
       (\rho-\rho_2)(\rho\rho_i-\rho_{i+1})}},
\\
      {\displaystyle
       {(\sigma^2 -\sigma_2)(\sigma_{i-1}-\sigma_i) \over
	(\sigma-\sigma_2)(\sigma_{i-1} - \sigma\sigma_i)}}
	&=& 
       {\displaystyle
       {(\rho^2 -\rho_2)(\rho_{i-1}-\rho_i) \over
	(\rho-\rho_2)(\rho_{i-1} - \rho\rho_i)}}. 
\end{eqnarray}
  \qed
\begin{th} \label{NT5.6}                             
Let $\Gamma = (X,R)$ denote a tight distance-regular graph with diameter
$d \ge 3$, and eigenvalues $\theta_0 > \theta_1 > \cdots >\theta_d$.
Let 
$\sigma_0,\sigma_1,\dots,\sigma_d$ denote the cosine sequence for
 $\theta_1$ or $\theta_d$.
 Fix adjacent vertices $x,y \in X$. Then with the notation of
Definition \ref{ZD4.1} we have the following (i), (ii).
\begin{description} \itemsep -1pt
\item{(i)} For all integers $i$ ($1\le i \le d-1$), and for all 
           $z \in D_{i}^{i}$,
           \begin{eqnarray} 
           |\Gamma(z) \cap D_{i-1}^{i-1}| \;&=&\; c_i \,
           {(\sigma_i-\sigma_{i+1})(\sigma \sigma_{i-1} -\sigma_i)
           \over 
           (\sigma_{i-1}-\sigma_i)(\sigma\sigma_i-\sigma_{i+1})},
                                                   \label{NT5.6B} 
\\
           |\Gamma(z) \cap D_{i+1}^{i+1}| \;&=&\; b_i \,
           {(\sigma_{i-1}-\sigma_i)(\sigma_i -\sigma\sigma_{i+1})
           \over 
           (\sigma_i-\sigma_{i+1})(\sigma_{i-1}-\sigma\sigma_i)}.
                                                     \label{NT5.6C}
           \end{eqnarray} 
\item{(ii)} For all integers $i$ ($2\le i \le d$), and for all 
           $z \in D_{i-1}^{i} \cup D_{i}^{i-1}$,
           \begin{equation}
          |\Gamma(z) \cap D_{i-1}^{i-1}| \;=\;
          a_{i-1}\,{{(1- \sigma) 
          (\sigma_{i-1}^2- \sigma_{i-2}\sigma_{i})}
          \over {(\sigma_{i-1} - \sigma_{i})
	         (\sigma_{i-2} - \sigma\sigma_{i-1})}}.    \label{NT5.6A}
          \end{equation}
\end{description} 
\end{th}

\proof 
(i) To prove (\ref{NT5.6B}), we assume $i\geq 2$; otherwise 
both sides are zero.
 Let $\alpha_i$ denote the expression on the right in
(\ref{NT5.5A}).
Let $N$ denote the number of ordered pairs
$uv$ such that 
$$
u \in  \Gamma_{i-1}(z) \cap D_{1}^{1}, \qquad 
           v\in \Gamma(z) \cap D_{i-1}^{i-1}, \qquad 
	   \partial(u,v)=i-2. 
$$
We compute $N$ in two ways. On one hand, by (\ref{NT5.5A}), there are precisely
$\alpha_i$ choices for $u$, and given $u$, there are precisely 
$c_{i-1}$ choices for $v$, so
\begin{equation}
N = \alpha_i c_{i-1}. \label{Nstep1}
\end{equation}
On the other hand, there are precisely $
\;|\Gamma(z) \cap D_{i-1}^{i-1}|\;$  choices for $v$, and given $v$, there
are precisely $\alpha_{i-1}$ choices for $u$, so
\begin{equation}
N = 
|\Gamma(z) \cap D_{i-1}^{i-1}|
\alpha_{i-1} . \label{Nstep2}
\end{equation}
Observe by Lemma \ref{L5.1}, Lemma \ref{AL3.9},  
and (\ref{NT5.5A})
that $\alpha_{i-1}\not=0$;
combining this with (\ref{Nstep1}), (\ref{Nstep2}), we find
$$
|\Gamma(z) \cap D_{i-1}^{i-1}| = c_{i-1}\alpha_i\alpha^{-1}_{i-1}.
$$
Eliminating $\alpha_{i-1}$, $\alpha_i$ in the above line  using 
(\ref{NT5.5A}), we obtain (\ref{NT5.6B}), as desired. 
\smallskip
\noindent 
Concerning  (\ref{NT5.6C}),  first assume $i=d-1$. We show
both sides
of (\ref{NT5.6C}) are zero. To see the  left side is zero, 
recall $a_d=0$  by Theorem \ref{ZZT4.9}, forcing $p^1_{dd}=0$ by
Lemma \ref{AL2.13}, so $D^d_d=\emptyset $ by the last line
in Definition \ref{ZD4.1}.
 The right side of (\ref{NT5.6C}) is zero since
the factor $\,\sigma_{d-1}-\sigma \sigma_d\,$ in the numerator is
zero by Lemma  \ref{TECH}(vi).
We now show 
 (\ref{NT5.6C})
for  $i\leq d-2$. Let $\beta_i$ denote the expression on the right
 in (\ref{NT5.5B}).
Let $N'$ denote the number of ordered pairs
$uv$ such that 
$$
u \in  \Gamma_{i+1}(z) \cap D_{1}^{1}, \qquad 
           v\in \Gamma(z) \cap D_{i+1}^{i+1}, \qquad 
	   \partial(u,v)=i+2.
$$
We compute $N'$ in two ways. On one hand, by (\ref{NT5.5B}), there are precisely
$\beta_i$ choices for $u$, and given $u$, there are precisely 
$b_{i+1}$ choices for $v$, so
\begin{equation}
N' = \beta_i b_{i+1}. \label{Npstep1}
\end{equation}
On the other hand, there are precisely $
\;|\Gamma(z) \cap D_{i+1}^{i+1}|\;$  choices for $v$, and given $v$, there
are precisely $\beta_{i+1}$ choices for $u$, so
\begin{equation}
N' = 
|\Gamma(z) \cap D_{i+1}^{i+1}|
\beta_{i+1} . \label{Npstep2}
\end{equation}
 Observe by Lemma \ref{L5.1}, Lemma  \ref{AL3.9}, and  
 (\ref{NT5.5B})
that $\beta_{i+1}\not=0$;
combining this with (\ref{Npstep1}), (\ref{Npstep2}), we find
$$
|\Gamma(z) \cap D_{i+1}^{i+1}| = b_{i+1}\beta_i\beta^{-1}_{i+1}.
$$
Eliminating $\beta_{i}$, $\beta_{i+1}$ in the above line  using 
(\ref{NT5.5B}), we obtain (\ref{NT5.6C}), as desired. 

\noindent (ii) Let $\gamma_i$ denote the expression on the right in 
(\ref{ZZT4.2A}), and let $\delta_i$ denote the expression on the right in 
(\ref{NT5.6B}). Let $N''$ denote the number of ordered pairs $uv$ such
that 
$$
u \in  \Gamma_{i-1}(z) \cap D_{1}^{1}, \qquad 
           v\in \Gamma(z) \cap D_{i-1}^{i-1}, \qquad 
	   \partial(u,v)=i-2.
$$
We compute $N''$ in two ways. On one hand, by Theorem \ref{ZZT4.2}(ii), there
are precisely $\gamma_i$ choices for $u$. Given $u$, we find by (\ref{NT5.6B})
(with $\,x $ and $i$ replaced by $u$ and $i-1$, respectively) that
there are precisely $c_{i-1}-\delta_{i-1}$ choices for $v$; consequently 
\begin{equation}
N''= \gamma_i(c_{i-1}-\delta_{i-1}).                         \label{Nppstep1}
\end{equation}
On the other hand, there are precisely  
$\vert\Gamma(z) \cap D_{i-1}^{i-1}\vert$ choices for $v$, and given $v$,
there are precisely $\alpha_{i-1}$ choices for $u$, where $\alpha_{i-1}$
is from the proof of (i) above. Hence
\begin{equation}
N''= \vert\Gamma(z) \cap D_{i-1}^{i-1}\vert \alpha_{i-1}.  \label{Nppstep2}
\end{equation}
Combining (\ref{Nppstep1}), (\ref{Nppstep2}), 
$$
\vert\Gamma(z) \cap D_{i-1}^{i-1}\vert =  
\gamma_i(c_{i-1}-\delta_{i-1})\alpha^{-1}_{i-1}.
$$
Eliminating $\alpha_{i-1}$, $\gamma_i$, $\delta_{i-1}$ in the above line
using (\ref{NT5.5A}), (\ref{ZZT4.2A}), (\ref{NT5.6B}),
respectively,  and simplifying the 
result using Theorem 
\ref{ZZT4.9}(ii), we obtain (\ref{NT5.6A}), as desired.
\qed

\begin{df}\label{ND4.4} \label{NND5.1} \label{ZD4.4} 
Let $\Gamma = (X, R)$ denote a distance-regular graph with diameter
$d \geq 3$, and fix adjacent vertices  $x,y \in X$.
\begin{description}
\item{(i)} For all integers $\;i,j$ we define the vector 
           $\;w_{ij}= w_{ij}(x,y)\;$ by 
           \begin{equation}
           w_{ij} = \sum_{z\in D_i^j} \hat z\;,
                                             \qquad \label{NND5.1B}
           \end{equation}
	   where $\;D^j_i=D^j_i(x,y)\;$ is from (\ref{ZD4.1A}). 
\item{(ii)} Let ${\cal L}$ denote the set of ordered pairs
           \begin{equation}
           {\cal L} \;=\; 
           \{ ij\ |\ 0\le i,j \le d,\ p_{ij}^1 \not= 0 \}.
                                             \qquad \label{NND5.1A}
           \end{equation}
           We observe that for all integers $i,j$, $w_{ij}\not=0$  if
	   and only if $\;ij\in {\cal L}$.
\item{(iii)} We define  the vector space  $W= W(x,y)$ by
	   \begin{equation} 
           W\;=\; \Span\{w_{ij}\; | \; ij\in {\cal L} \}.
           \label{NND5.1C}
           \end{equation}
\end{description}
\end{df}

\begin{lemma} \label{NL4.7}  \label{NNL5.2} \label{ZL4.5} 
Let $\;\Gamma=(X,R)\;$ denote a distance-regular graph with diameter 
$\;d\geq 3$, and assume $a_1 \not= 0$. Then 
\begin{description}
\item{(i)} \ \ ${\cal L} = \{i-1,i \;|\; 1 \le i \le d \}
                       \cup \{i,i-1 \;|\; 1 \le i \le d \}
                       \cup \{ii \;|\; 1 \le i \le e \}$, \\
            where $e = d-1$ if $a_d=0$ and $e=d$ if $a_d \not= 0$.
\item{(ii)} 
            \begin{equation}
	       |{\cal L}| = 
            \left\{ \begin{array}{ll}
                    3d  & \mbox{if\ \  $a_d \not= 0$,} \\
                    3d-1& \mbox{if\ \  $ a_d = 0$.}
                    \end{array} \right.
                                             \qquad \label{NNL5.2A}
            \end{equation}
\item{(iii)} Let $x,y$ denote adjacent vertices in $X$, and let
	  $W=W(x,y)$ be as in (\ref{NND5.1C}).  Then  
            \begin{equation}
            \dim W \;=\; 
            \left\{ \begin{array}{ll}
                    3d  & \mbox{if\ \  $a_d \not= 0$,} \\
                    3d-1& \mbox{if\ \  $ a_d = 0$.}
                    \end{array} \right.
                                             \qquad \label{NNL5.2B}
            \end{equation}
\end{description}
\end{lemma}

\proof Routine application  of Lemma \ref{ZZL2.11} and Lemma \ref{AL2.13}.
 \qed

\begin{lemma}  \label{NNL6.2}                        
Let $\Gamma = (X,R)$ denote a distance-regular graph with diameter
$d \ge 3$,  fix adjacent vertices $x, y \in X$,
and let the vector space
$W=W(x,y)$ be as in (\ref{NND5.1C}). Then 
the following are equivalent.
\begin{description} \itemsep -3pt
\item{(i)} The vector space $W$ 
            is $A$-invariant.
\item{(ii)} For all integers $i,j,r,s$ ($ij \in {\cal L}$ and 
            $rs \in {\cal L}$), and for all $z \in D_{i}^{j}$, 
            the scalar $|\Gamma(z) \cap D_{r}^{s}|$ is a constant
            independent of $z$. 
%
\item{(iii)} The following conditions hold.
            \begin{description} \itemsep -3pt
            \item{(a)} For all integers $i$ ($1 \le i \le d$), 
                       and for all $z \in D_{i}^{i}$, the scalars
                       $|\Gamma(z) \cap D_{i-1}^{i-1}|$ and 
                       $|\Gamma(z) \cap D_{i+1}^{i+1}|$ 
                       are constants independent of $z$.
            \item{(b)} For all integers $i$ ($2 \le i \le d$), 
                       and for all $z \in D_{i-1}^{i}\cup D^{i-1}_{i}$,
		       the scalar $|\Gamma(z)\cap D_{i-1}^{i-1}|$ 
                       is a constant independent of $z$.
	 \end{description}
\end{description}
\end{lemma}
\proof (i) \iff (ii) Routine. 
\\
(ii) \implies (iii) Clear.
\\
(iii) \implies (ii) Follows directly from Lemma \ref{NL4.6}. 
\qed


\begin{df}                                        
\label{homog}
Let $\Gamma = (X,R)$ denote a distance-regular graph with diameter
$d \ge 3$.
For each edge  $xy \in R$, 
the graph $\Gamma$ is said to be {\bf 1-homogeneous with respect to 
$xy$} whenever (i)--(iii) hold in Lemma \ref{NNL6.2}.
The graph $\Gamma$ is said to be {\bf 1-homogeneous} whenever it is 
1-homogeneous with respect to all edges in $R$.
%
\end{df}

\noindent 
 
\begin{th} \label{NT4.11}                         
Let $\Gamma = (X,R)$ denote a distance-regular graph
with diameter
$d \ge 3$.   Then the following are equivalent.
\begin{description} \itemsep -3pt
\item{(i)}   $\;\;\Gamma$ is tight,
\item{(ii)} $\;a_1\not=0$,  $a_d=0$, and  $\Gamma $ is 1-homogeneous,
\item{(iii)} $a_1\not=0$, $a_d=0$, and $\Gamma$ is 1-homogeneous with respect
to at least one edge.
\end{description}
\end{th}

\proof (i) \implies (ii) 
 Observe $a_1\not=0$ by Proposition \ref{ZP3.7}, and 
 $a_d=0$  by Theorem \ref{ZZT4.9}. 
 Pick any   
 edge $xy \in R$.  By Theorem \ref{NT5.6}, we find   
   conditions (iii)(a), (iii)(b) hold in 
 Lemma \ref{NNL6.2}, so $\Gamma$ is 
 1-homogeneous with respect to $xy$ by Definition
\ref{homog}. Apparently $\Gamma$ is 1-homogeneous with respect to every
edge, so $\Gamma $ is 1-homogeneous. 
%
%
\\
(ii) \implies (iii) Clear.
\\
(iii) \implies (i) Suppose $\Gamma $ is 1-homogeneous with respect
to the edge $xy \in R$.  We show $xy $ is tight with respect to both
$\theta_1, \theta_d$. To do this, we show the tightness
$t=t(x,y)$ from Definition \ref{deftight}
 equals 2.
Consider the vector space $W=W(x,y)$ from 
(\ref{NND5.1C}), and
the vector space $H$ from 
 (\ref{NL4.2B}). 
Observe  $W$  is    
$A$-invariant by Lemma \ref{NNL6.2}, and $W$  contains $H$, so
it contains $MH$, where $M$ denotes the Bose-Mesner algebra of $\Gamma $.
The space $W$ has dimension $3d-1$ by
(\ref{NNL5.2B}), so
 $MH$ has dimension at most $3d-1$. Applying  
 (\ref{NL4.2A}), we find $t\geq  2$. From the discussion at the end
 of Definition \ref{deftight}, we observe $t=2$, and that
$xy$ is tight with respect to both  
$\theta_1$, $\theta_d$.
Now  $\Gamma$ is tight in view of Corollary \ref{ZC3.6}(iv)  and
Definition \ref{AD3.8}.   \qed

\section{The local graph}                             

\begin{definition}                                 \label{localgraph1}
Let $\Gamma=(X,R)$ denote a distance-regular graph with diameter 
$d\ge 3$. For each vertex $x\in X$, we let $\Delta = \Delta(x)$ denote 
the vertex subgraph of $\Gamma$ induced on $\Gamma (x)$. 
We refer to $\Delta $ as the {\bf local graph} associated with $x$. 
We observe $\Delta $ has $k$ vertices, and is regular with valency 
$a_1$. We  further observe $\Delta $ is not a clique.
\end{definition}

\noindent 
In this section, we show the local graphs of tight distance-regular
graphs are strongly-regular. We begin by  recalling the definition
and some basic properties of strongly-regular graphs.

\begin{definition}\cite[p.3]{BCN}
\label{SRDEF}
A graph $\Delta $ is said to be {\bf strongly-regular} with {\rm parameters}
$(\nu,\kappa,\lambda, \mu)$ whenever $\Delta $ has $\nu$ vertices and is
regular with valency $\kappa$, adjacent vertices of $\Delta $ have precisely 
$\lambda $ common neighbors, and distinct  non-adjacent vertices of $\Delta $
have precisely $\,\mu\,$ common neighbors.
\end{definition}

\begin{lemma}\cite[Thm. 1.3.1]{BCN}
\label{SRGFORM}
Let $\,\Delta \,$ denote a connected  strongly-regular graph with
parameters 
$\,(\nu,\kappa,\lambda, \mu)$, and assume $\Delta$  is not a clique. Then
$\,\Delta \,$ has precisely three distinct eigenvalues, one of which
is $\kappa$. Denoting the others by 
 $\,r, s$,
\begin{equation}
\nu = {(\kappa-r)(\kappa-s)\over \kappa+rs},\ \ \ \ \ \ \ \lambda = \kappa + r + s + rs, 
\ \ \ \ \ \ \ \mu = \kappa + r s.
\label{SRGFORM1}
\end{equation}
\noindent
The multiplicity of $\kappa $  as an eigenvalue of $\Delta $
 equals   1. 
The multiplicities with which  $\,r, s\,$ appear as eigenvalues of $\,\Delta\, $
are given by 
\begin{equation}
mult_{r} = {\kappa (s+1)(\kappa-s) \over \mu (s - r)},
\qquad \qquad mult_{s} = {\kappa (r+1)(\kappa-r) \over \mu (r - s)}.
\label{SRGFORM2}
\end{equation}
\end{lemma}

%
%
%

%

\begin{th}\label{DELTASR}                       
Let $\Gamma = (X,R) $ denote a tight distance-regular graph with 
diameter $d\ge 3$, and eigenvalues $\theta_0>\theta_1>\cdots >\theta_d$.
Pick $\theta\in \lbrace \theta_1, \theta_d\rbrace$, let $\sigma$, 
$\sigma_2$ denote the first and second cosines for $\theta $, 
respectively, and let $\eps$ denote the associated auxiliary parameter 
from (\ref{ZZT4.4A}). Then for any vertex $\;x\in X$, the local graph 
$\;\Delta=\Delta(x)\;$ satisfies (i)--(iv)  below. 
\begin{description} \itemsep -3pt
\item{(i)} \ \ $\Delta$ is strongly-regular with  parameters 
               $\;(k,a_1,\lambda,\mu)$, where $k$ is the valency of 
               $\,\Gamma$, and    
\begin{eqnarray} 
a_1 &=&  
-{(1-\sigma_2)(1+\sigma)(1-\eps)\over (\sigma-\sigma_2)(1-\eps \sigma)},
\label{SRGA1} \\
\lambda &=&  a_1{{2\sigma}\over {1+\sigma}}\;-\; 
a_1{1-\sigma\over 1+\sigma}\,{\sigma_2\over \sigma-\sigma_2}\;-\; 
{1-\sigma_2\over \sigma-\sigma_2}, 
\label{LAMBA} \\
\mu &=& {a_1 \over 1+\sigma}\,{\sigma^2-\sigma_2\over \sigma-\sigma_2}.
\label{MUA}
\end{eqnarray}
\item{(ii)} \ \  $\Delta $ is connected and not a clique. 
\item{(iii)} \ \ The distinct eigenvalues of $\,\Delta \,$ are 
                 $\,a_1$, $r$, $s$, where 
\begin{eqnarray}                                            \label{MUC}
r= {{a_1 \sigma }\over {1+\sigma}}, \qquad \qquad 
s= -{{1-\sigma_2}\over {\sigma-\sigma_2}}.  
\end{eqnarray}
\item{(iv)} The multiplicities of $\,r$, $s\,$ are given by
\begin{eqnarray}
mult_r=  {{(1+\sigma)(\sigma-\eps)}\over {\sigma_2-\sigma^2}},
\qquad \qquad 
mult_s=  -{(1-\eps)(1+\sigma)(\sigma_2-\eps\sigma)\over 
           (\sigma_2-\sigma^2)(1-\eps \sigma)}. \qquad 
\label{SRGLOCM}
\end{eqnarray}
\end{description}
\end{th}

\proof 
(i) Clearly $\,\Delta\,$ has $k$ vertices and is regular with valency 
$a_1$. The formula (\ref{SRGA1}) is from Theorem \ref{ZZT4.9}(ii). 
Pick distinct vertices $\,y,\, z\in \Delta$. 
We count the number of common neighbors of $\,y,\,z\,$ in $\Delta$.  
First suppose $\,y,\,z\,$ are adjacent. By (\ref{NT4.12B}) (with $i=1$)
we find $\,y,\,z\,$ have precisely $\,\lambda \,$ common neighbors in 
$\,\Delta$, where $\,\lambda\, $ is given in (\ref{LAMBA}).
Next  suppose $\,y,\,z\,$ are not adjacent. By (\ref{ZZT4.2A}) 
(with $i=2$), we find $\,y,\,z\,$ have precisely $\,\mu\,$ common 
neighbors in $\,\Delta$, where $\,\mu\,$ is given in (\ref{MUA}). 
The result now follows in view of Definition \ref{SRDEF}.

\noindent 
(ii)  
We saw in Definition \ref{localgraph1} that $\Delta$ is not a clique. 
Observe  the scalar $\,\mu\,$ in (\ref{MUA}) is not zero, since 
$\,a_1\ne 0\,$ by Proposition \ref{ZP3.7}, and since 
$\,\sigma^2\ne\sigma_2\,$ by Lemma \ref{AL3.9}(ii),(iii). 
It follows $\,\Delta\, $ is connected.

\noindent (iii) The scalar $a_1$ is an eigenvalue of $\Delta $ by 
Lemma \ref{SRGFORM}. Using (\ref{LAMBA}), (\ref{MUA}), we find  
the scalars  $r,\,s$ in   (\ref{MUC}) satisfy  
$$
\lambda = a_1 +r + s + rs, \qquad \qquad \mu= a_1 + rs. 
$$
Comparing this with the two equations on the right in (\ref{SRGFORM1}),
we find the scalars $r, s$ in (\ref{MUC}) are the remaining eigenvalues
of $\Delta $.

\noindent 
(iv) By (\ref{SRGFORM2}) and (i) above,
$$
mult_{r} = {a_1(s+1)(a_1-s) \over \mu (s - r)},
\qquad \qquad mult_{s} = {a_1(r+1)(a_1-r) \over \mu (r - s)}.
$$ 
Eliminating $\,a_1,\, \mu, \, r,\, s \,$ in the above equations using 
(\ref{SRGA1}), (\ref{MUA}), (\ref{MUC}), we routinely obtain 
(\ref{SRGLOCM}).  \qed

\begin{definition} \label{BPM}          
Let $\Gamma $ denote a distance-regular graph with diameter $\,d\geq 3$,
and  eigenvalues $\theta_0> \theta_1 >\cdots > \theta_d$. 
We define 
$$
b^- \;:=\; -1\;-\;{{b_1}\over {1+\theta_1}}, \qquad \qquad 
b^+ \;:=\; -1\;-\;{{b_1}\over {1+\theta_d}}.
$$
We recall $a_1-k\leq \theta_d < -1<\theta_1$ by Lemma \ref{ZL2.0},  
so $\;b^-<-1$, $\;b^+\geq 0$.
\end{definition}

\begin{th}\label{ZT5.8} \label{AT6.2}                  
Let $\Gamma=(X,R)$ denote a distance-regular graph with 
diameter $d \geq 3$.
 Then the following are equivalent. 
\begin{description} \itemsep -3pt
\item{(i)} \ \ $\Gamma$ is tight.
\item{(ii)} \ For all $x \in X$, the local graph $\Delta(x)$ is 
          connected  strongly-regular with 
	  eigenvalues $a_1$, $b^+$, $b^-$. 
\item{(iii)} There exists  $x\in X $ for which the local graph 
	     $\Delta(x)$ is connected strongly-regular with   
	     eigenvalues $a_1$,  
	     $b^+$,  $b^-$.
\end{description}
\end{th}

\proof  
\noindent (i) \implies (ii)  
Pick any $x\in X$, and let $\Delta=\Delta(x)$ denote the
local graph. By Theorem \ref{DELTASR},  
the graph  $\,\Delta\,$ is connected and strongly-regular. The eigenvalues 
 of $\,\Delta\, $ other than $\,a_1\,$ are given by (\ref{MUC}),
 where for convenience  we take  the eigenvalue $\,\theta\, $ involved
 to be $\,\theta_1$. 
Eliminating $\,\sigma$, $\sigma_2 $ in (\ref{MUC})  using
$\theta_1 = k\sigma $ and Lemma \ref{TECH}(i), 
and simplifying the results using equality
in the fundamental bound (\ref{ZT3.5A}),  we routinely find $r=b^+$, $s=b^-$. 

\noindent (ii) \implies (iii) Clear.

\noindent (iii) \implies (i) Since $\Delta = \Delta(x)$ is connected, its
valency $a_1$ is not zero. In particular $\Gamma$ is not bipartite. 
The graph $\Delta $ is not a clique, so 
 (\ref{SRGFORM1}) holds for $\Delta$. Applying the equation on the left
 in that line, we obtain
\begin{eqnarray}
k(a_1 +b^+b^-) \;=\; (a_1-b^+)(a_1- b^-).
\label{FBVER2}
\end{eqnarray}
Eliminating $\,b^+$, $\,b^-\,$ in (\ref{FBVER2}) using 
Definition \ref{BPM}, and simplifying the result, 
we routinely obtain equality in the fundamental bound (\ref{ZT3.5A}).  
Now $\Gamma $ is tight, as desired. \qed

%
%
%
%

\section{Examples of tight distance-regular graphs} 

The following examples (i)-(xii) are tight distance-regular graphs with 
diameter at least 3.
In each case we give
the intersection array, the second largest eigenvalue $\theta_1$,
and the least eigenvalue $\theta_d$, together with 
their respective cosine sequences $\{\sigma_i\}$, $\{\rho_i\}$, and
the 
 auxiliary parameter $\eps$ for $\theta_1$.
%
%
Also, we give  the parameters and nontrivial eigenvalues of 
the local graphs.

\bigskip\noindent
{\bf (i)} The {\bf Johnson graph $J(2d,d)$} has diameter $d$ and 
intersection numbers $a_i=2i(d-i)$, $b_i=(d-i)^2$, $c_i = i^2$ for 
$i=0,\dots,d$, cf. \cite[p. 255]{BCN}. It is distance-transitive, 
an antipodal double-cover,  and $Q$-polynomial with respect to $\theta_1$.
%

\noindent Each local graph is a {\bf lattice graph $K_d\times K_d$},
%
%
 with parameters $(d^2,2(d-1),d-2,2)$
and nontrivial eigenvalues $\,r=d-2$, $\,s=-2$, cf. \cite[p. 256]{BCN}.

\medskip\noindent
{\bf (ii)} The {\bf halved cube $\frac12 H(2d,2)$}  has diameter $d$ and
intersection numbers $a_i=4i(d-i)$, $b_i=(d-i)(2d-2i-1)$, $c_i=i(2i-1)$ 
for $i=0,\dots d$, cf. \cite[p. 264]{BCN}. It is distance-transitive,
an antipodal double-cover, and $Q$-polynomial with respect to $\theta_1$.
%

\noindent Each local graph is a  
 Johnson graph $J(2d,2)$, with parameters 
$(d(2d-1),4(d-1),2(d-1),4)$ and nontrivial eigenvalues 
$\,r=2d-4$, $\,s=-2$, cf. \cite[p. 267]{BCN}.


\medskip\noindent
{\bf (iii)} The {\bf Taylor graphs} are nonbipartite double-covers of complete 
graphs,  i.e., distance-regular graphs with intersection array of the form
$\{k,c_2,1;1,c_2,k\}$, where $c_2<k-1$.
They have diameter 3, and are  
 $Q$-polynomial with respect to
both $\theta_1$, $\theta_d$. These eigenvalues are given by
 $\theta_1=\alpha$, 
$\theta_d=\beta$, where  
$$
\alpha + \beta =k-2c_2-1, \qquad \qquad \alpha \beta = -k,
$$ 
and $\alpha > \beta$.  
See Taylor \cite{Ta}, and Seidel and Taylor \cite{ST} for more details. 
 
\noindent Each local graph  is strongly-regular
with  parameters $(k,a_1,\lambda,\mu)$, where $\;a_1=k-c_2-1$, 
$\;\lambda = (3a_1-k-1)/2$
and $\;\mu=a_1/2$. We note both $\;a_1, c_2 \,$ are even and $k$ is odd. 
The nontrivial  
eigenvalues of the local graph are 
$$
r = {{\alpha-1} \over 2},\ \ \ \ \ \ \ \ \ \ \  
s = {{\beta -1} \over 2}.
$$
%

%
%


\medskip\noindent 
{\bf (iv)} The graph {\bf $3.Sym(7)$} has intersection array 
$\{ 10,6,4,1;1,2,6,10\}$ and can be obtained from 
a sporadic Fisher group, 
%
%
cf.  \cite[pp. 397-400]{BCN}. It is sometimes called the 
Conway-Smith graph. It is distance-transitive, an antipodal 
3-fold cover, and is not $Q$-polynomial.   
%

\noindent Each local graph is a  {\bf Petersen graph}, with parameters
$(10,3,0,1)$ and nontrivial eigenvalues $\,r = 1$, $\,s=-2$,
see \cite{Hal1}, \cite[13.2.B]{BCN}.

\medskip\noindent
{\bf (v)} The graph {\bf $3.O_6^-(3)$} has intersection array 
$\{45,32,12,1;1,6,32,45\}$ and can be obtained from 
a sporadic Fisher group, 
%
%
cf.  \cite[pp. 397-400]{BCN}. It is distance-transitive, an antipodal 
3-fold cover,  and is not $Q$-polynomial. 
%

\noindent 
Each local graph is a  {\bf generalized quadrangle $GQ(4,2)$},
 with parameters $(45,12,3,3)$ and  nontrivial
eigenvalues $\,r=3$,  $\,s=-3$.
See \cite[p. 399]{BCN}.

\medskip\noindent 
{\bf (vi)} The graph {\bf $3.O_7(3)$} has intersection array 
$\{ 117,80,24,1;1,12,80,117\}$ and can be obtained from 
a sporadic Fisher group,
%
%
cf. \cite[pp. 397-400]{BCN}. It is distance-transitive, an antipodal 
3-fold cover,  and is not $Q$-polynomial. 
%

\noindent Each local graph is strongly-regular  with parameters 
$(117,36,15,9)$,
 and nontrivial eigenvalues
$\,r=9$, $\,s=-3$.
\cite[13.2.D]{BCN}.

\medskip\noindent 
{\bf (vii)} The graph {\bf $3.Fi_{24}$} has intersection array
$\{31671,28160,2160,1;1,1080,28160,31671\}$ and can be obtained from
a sporadic Fisher group, 
%
%
cf.  \cite[pp. 397]{BCN}. It is distance-transitive, 
an antipodal 3-fold cover,  and is not $Q$-polynomial. 
%
%

\noindent 
Each local graph is strongly-regular  with parameters
$(31671,3510,693,351)$ and  nontrivial eigenvalues 
$\,r = 351$,  $\,s=-9$. They are related to $Fi_{23}$.


\medskip\noindent 
{\bf (viii)} The {\bf Soicher1 graph} has intersection array 
$\{ 56,45,16,1;1,8,45,56\}$, cf. \cite{Br3}, \cite[$\!$11.4I]{BCN2},
\cite{Soi}. It 
is distance-transitive, 
 an  antipodal 3-fold cover,  and is not $Q$-polynomial. 
%

\noindent Each local graph is a {\bf Gewirtz graph} with parameters
$(56,10,0,2)$ and nontrivial eigenvalues $\,r = 2$,  $\,s=-4$,
\cite[p.372]{BCN}.

\medskip\noindent 
{\bf (ix)} The {\bf Soicher2 graph} has intersection array 
$\{ 416,315,64,1;1,32,315,416\}$, cf. \cite{Soi} 
\cite[$\!$13.8A]{BCN2}.
It is distance-transitive,
an antipodal 3-fold cover,  and is not $Q$-polynomial. 
%

\noindent Each local graph is  strongly-regular with parameters 
$(416,100,36,20)$ and  nontrivial eigenvalues $\,r = 20$,  
$s=-4$. 
       

\medskip\noindent 
{\bf (x)} The {\bf Meixner1 graph} has intersection array 
$\{ 176,135,24,1;1,24,135,176\}$, cf. \cite{Mei} 
\cite[$\!$12.4A]{BCN2}.
It is distance-transitive,
an antipodal 2-fold cover,  and is $Q$-polynomial. 

\noindent Each local graph is  strongly-regular with parameters 
$(176,40,12,8)$ and  nontrivial eigenvalues $\,r = 8$,  
$s=-4$.

\medskip\noindent 
{\bf (xi)} The {\bf Meixner2 graph} has intersection array 
$\{ 176,135,36,1;1,12,135,176\}$, cf. \cite{Mei} 
\cite[$\!$12.4A]{BCN2}.
It is distance-transitive,
an antipodal 4-fold cover,  and is not $Q$-polynomial. 

\noindent Each local graph is  strongly-regular with parameters 
$(176,40,12,8)$ and  nontrivial eigenvalues $\,r = 8$,  
$s=-4$.


\medskip\noindent
{\bf (xii)} The {\bf Patterson graph} has intersection array
$\{ 280, 243,144,10; 1,8,90,280\}$, and can be constructed from the 
Suzuki group, see \cite[13.7]{BCN}. 
It is primitive and distance-transitive, but not $Q$-polynomial. 

\noindent
Each local graph is a  generalized quadrangle $GQ(9,3)$ with parameters 
$(280,36,8,4)$ and nontrivial eigenvalues $\,r = 8$, $\,s=-4$,
\cite[Thm. 13.7.1]{BCN}. 

{\small
\def\tablerule{\noalign{\hrule}}
\hskip -1.3cm
{
$$
\vbox{\offinterlineskip \tabskip=0pt
\halign{\vrule# \tabskip=3pt      
&\ns\hfil # \hfil \ns\ns\ns          
&\vrule#                          
&\hfil $#$ \hfil    \ns           
&\hfil $#$ \hfil    \ns           
&\vrule#                          
&\hfil $#$ \hfil                  
&\hfil $#$ \hfil                  
&\hfil $#$ \hfil                  
&\vrule#\tabskip=0pt\cr         
%
\tablerule height 15pt depth 10pt 
& Name && \theta_1 & \theta_d && \{\sigma_i\} & \{\rho_i\}  & \eps &\cr
\tablerule height 15pt depth 10pt 
& 
$J(2d,d)$ && d(d-2) & -d
&& 
  \sigma_i = {\displaystyle{d-2i \over d} }
& \rho_i   = {\displaystyle{(-\!1)^i\!\cdot\! 1\!\cdot \!2 
         \cdots i\over 
               d(d\!-\!1)\cdots (d\!-\!i\!+\!1)}}
& {\displaystyle {d+2\over d} }
&\cr 
\tablerule height 15pt depth 10pt 
%
& 
${1 \over 2} H(2d,2)$ && (2d\!-\!1)(d\!-\!2) &  -d
&& 
  \sigma_i = {\displaystyle{d\!-\!2i\over d} }
& \!\!\!\! 
  \rho_i   = {\displaystyle{(-\!1)^i\! \cdot\! 1\!\cdot\! 3 \cdots (2i\!-\!1)\ 
             \over(2d\!-\!1)(2d\!-\!3)\cdots(2d\!-\!2i\!+\!1)}}
& {\displaystyle {d+1 \over d-1} }
&\cr 
\tablerule height 15pt depth 10pt 
%
& 
Taylor && \alpha & \beta 
%
%
&& 
  {\displaystyle (1,{\alpha \over k}, {-\alpha \over k}, -1) }
& {\displaystyle (1,{\beta \over k}, {-\beta \over k}, -1) }
& {\displaystyle {k+1 \over \alpha \!-\! \beta} }
&\cr 
\tablerule height 15pt depth 10pt 
%
& 
3.Sym(7) && 5 &  -4 
&& 
  {\displaystyle (1,{1 \over 2}, 0, {-1 \over 4}, {-1 \over 2})   }
& {\displaystyle (1,{-2 \over 5}, {3 \over 10}, {-2 \over 5}, 1)  }
& {\displaystyle {4 \over 3} }
&\cr 
\tablerule height 15pt depth 10pt 
%
& 
$3.O_6^-(3)$ && 15 & -9
&& 
  {\displaystyle (1,{1 \over 3}, 0, {-1 \over 6}, {-1 \over 2}) }
& {\displaystyle (1,{-1 \over 5}, {1 \over 10}, {-1 \over 5}, 1) }
& 2
&\cr 
\tablerule height 15pt depth 10pt 
%
& 
$3.O_7(3)$ && 39 & -9
&& 
  {\displaystyle (1,{1 \over 3}, 0, {-1 \over 6}, {-1 \over 2}) }
& {\displaystyle (1,{-1 \over 13}, {2 \over 65}, {-1 \over 13}, 1) }
& {\displaystyle {5 \over 2} }
&\cr 
\tablerule height 15pt depth 10pt 
%
& 
$3.Fi_{24}$ && 3519 & -81
&& 
  {\displaystyle (1,{1 \over 9}, 0, {-1 \over 18}, {-1 \over 2}) }
& {\displaystyle (1,{-1 \over 391}, {5 \over 17204}, {-1 \over 391}, 1) }
& {\displaystyle {44 \over 5} }
&\cr 
\tablerule height 15pt depth 10pt 
%
& 
Soicher1 && 14 & -16
&& 
  {\displaystyle (1,{1 \over 4}, 0, {-1 \over 8}, {-1 \over 2}) }
& {\displaystyle (1,{-2 \over 7}, {1 \over 7}, {-2 \over 7}, 1) }
& 2
&\cr 
\tablerule height 15pt depth 10pt 
%
& 
Soicher2 && 104 & -16
&& 
  {\displaystyle (1,{1 \over 4}, 0, {-1 \over 8}, {-1 \over 2}) }
& {\displaystyle (1,{-1 \over 26}, {1 \over 91}, {-1 \over 26}, 1) }
& {\displaystyle {7 \over 2} }
&\cr 
\tablerule height 15pt depth 10pt 
%
& 
Meixner1 &&  44 & -16
&& 
  {\displaystyle (1,{1 \over 4}, 0, {-1 \over 4}, -1) }
& {\displaystyle (1,{-1 \over 11}, {1 \over 33}, {-1 \over 11}, 1) }
& {\displaystyle 3}
&\cr 
\tablerule height 15pt depth 10pt 
%
& 
Meixner2 &&  44 & -16
&& 
  {\displaystyle (1,{1 \over 4}, 0, {-1 \over 12}, {-1\over 3}) }
& {\displaystyle (1,{-1 \over 11}, {1 \over 33}, {-1 \over 11}, 1) }
& {\displaystyle 3}
&\cr 
\tablerule height 15pt depth 10pt 
%
& 
Patterson && 80 & -28 
&& 
  {\displaystyle (1,{2 \over 7}, {1 \over 21}, {-2 \over 63}, {-1 \over 9}) }
& {\displaystyle (1,{-1 \over 10}, {1 \over 45}, {-1 \over 54}, {5 \over 27}) }
& {\displaystyle  {8 \over 3} }
&\cr 
%
\tablerule
}}
$$
} 
} 

\bigskip \noindent
ACKNOWLEDGEMENT: We would like to thank Prof. Yoshiara for
mentioning that the Patterson graph satisfies the Fundamental Bound.




\baselineskip=\normalbaselineskip
\newpage
\section{Appendix A: 
1-homogeneous partitions of the known examples of the AT4 family
and the Patterson graph}
%
%
%
\def\G{\Gamma} {\rm  
In \cite{JK3} a tight nonbipartite antipodal distance-regular graph $\G$
with diameter four was parametrized by the eigenvalues $r$ and $-s$
of the local graphs and the size $t$ of its antipodal classes. The graph
$\G$ was called an {\de antipodal tight graph} of diameter four and 
with parameters $(r,s,t)$, and denoted by AT4$(r,s,t)$.}

\medskip\medskip
\centerline{\psfig{figure=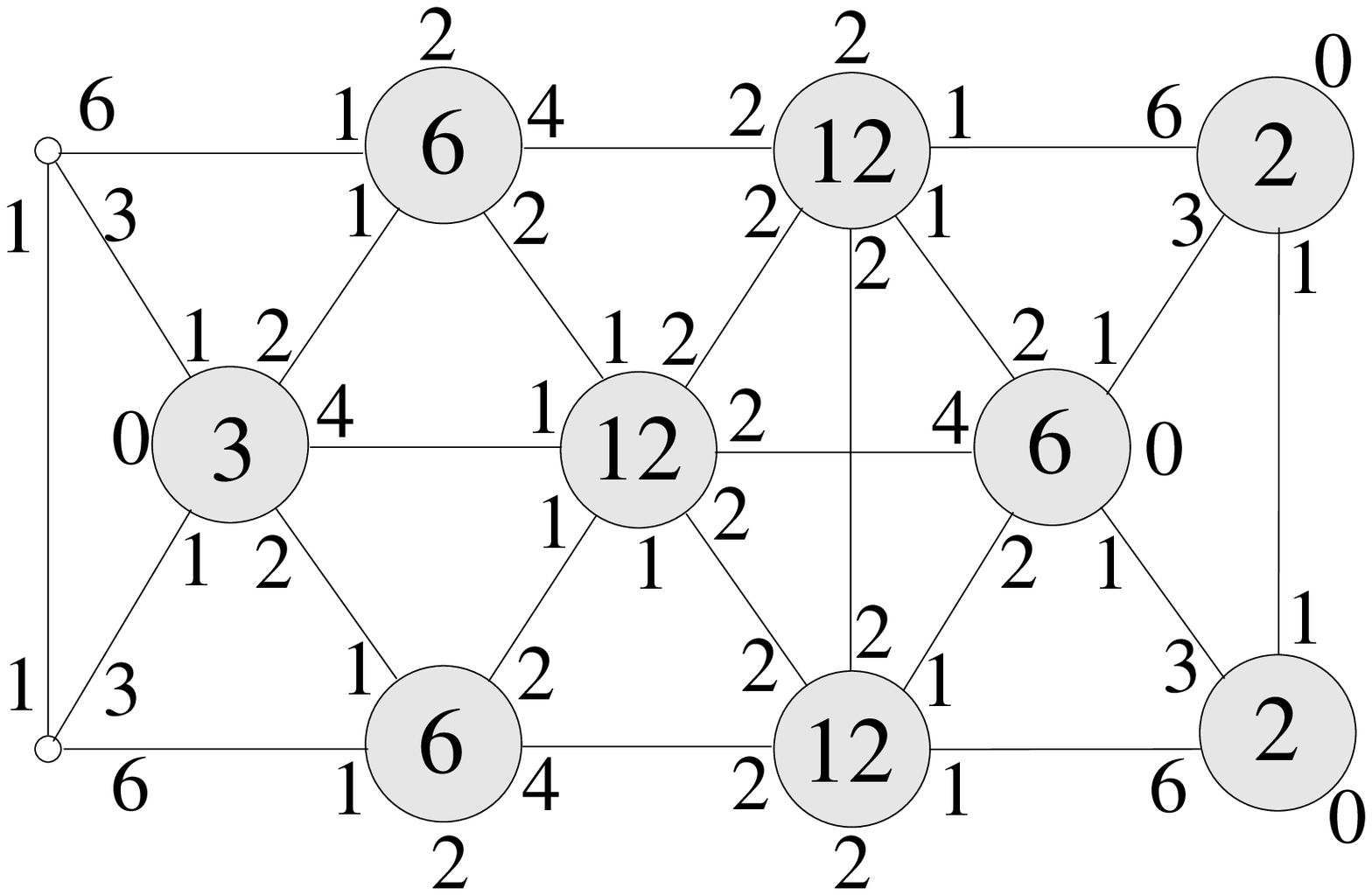,height=28mm}\ \ \ 
            \psfig{figure=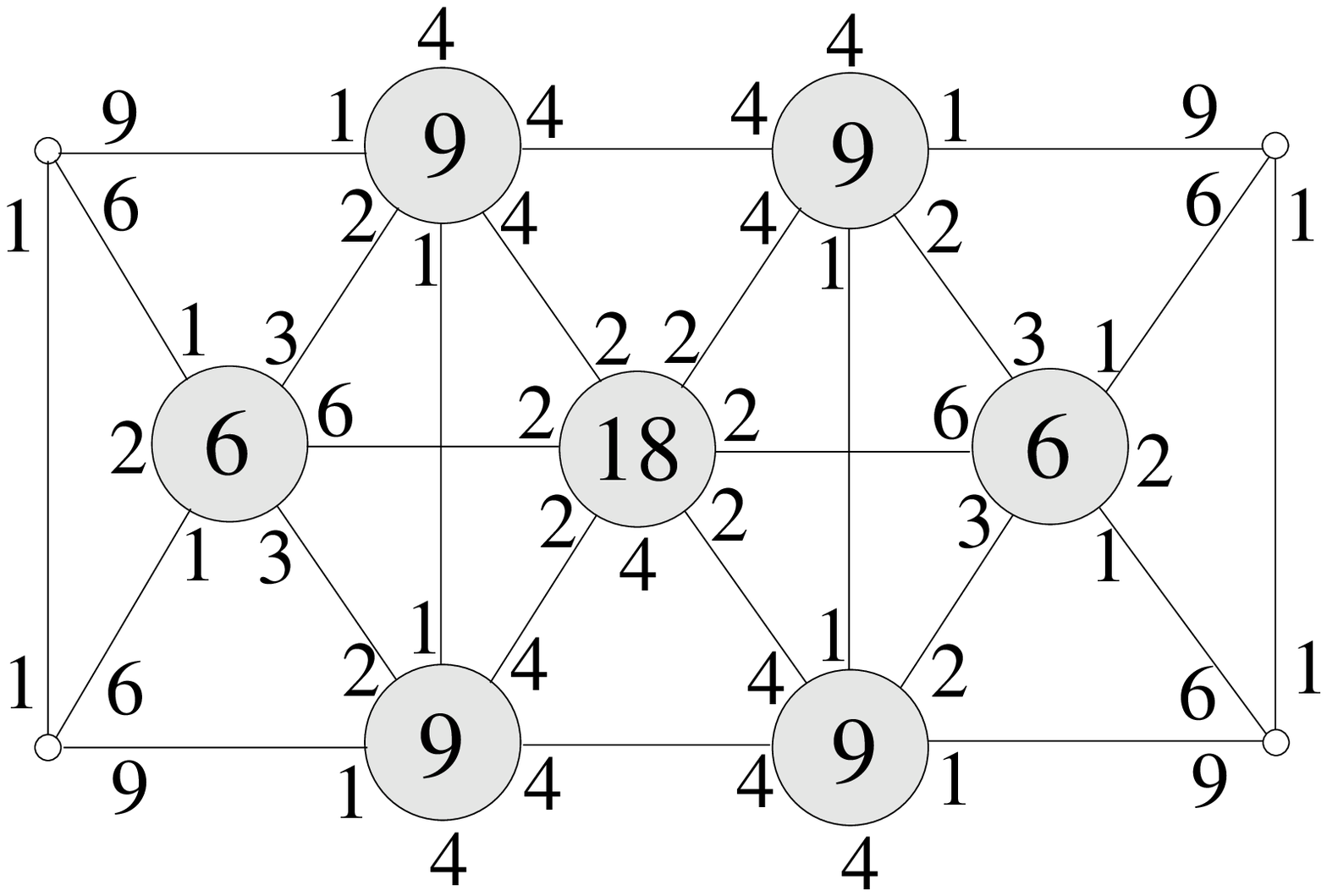,height=28mm}\ \ \ 
            \psfig{figure=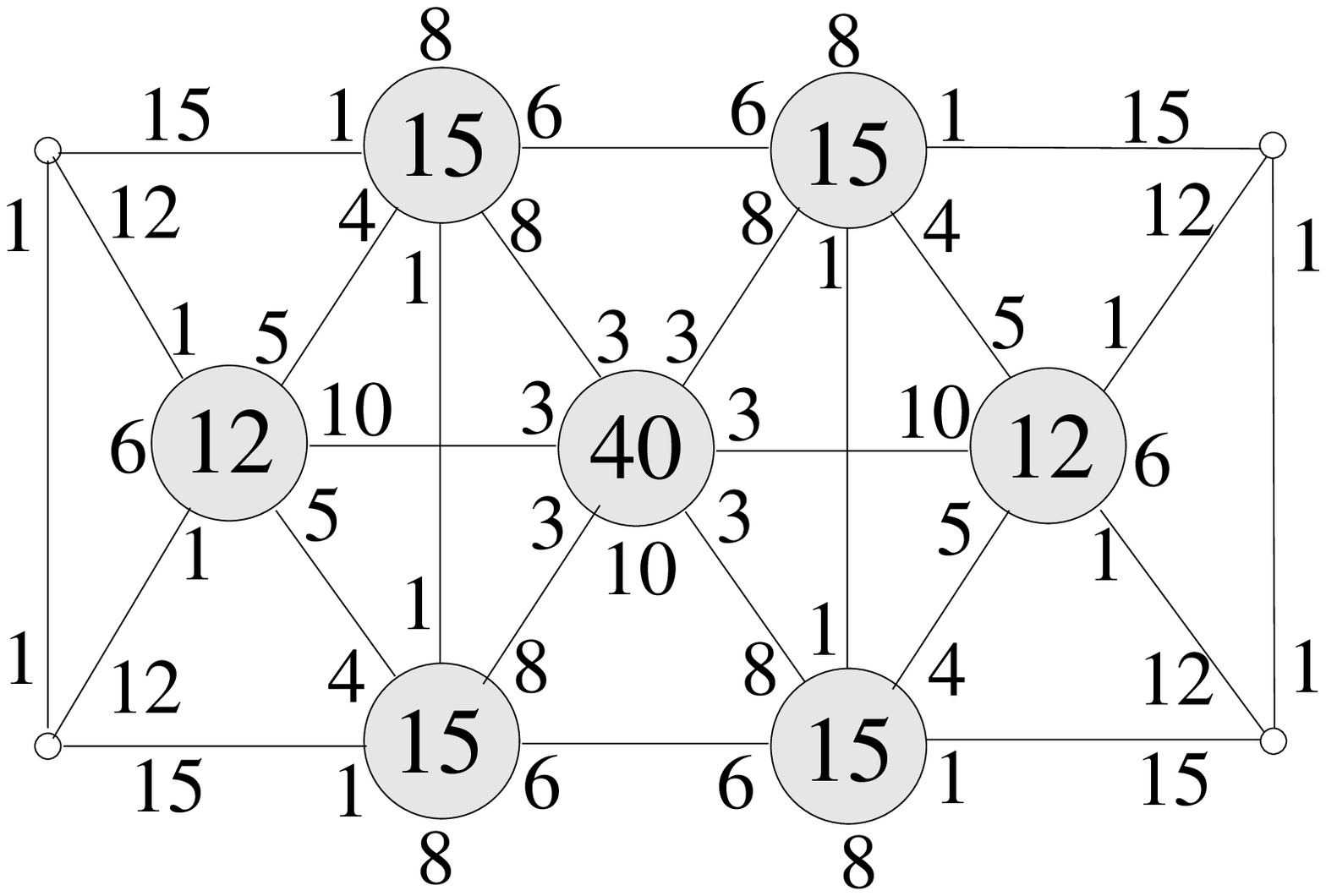,height=28mm}\ \ \ 
            \psfig{figure=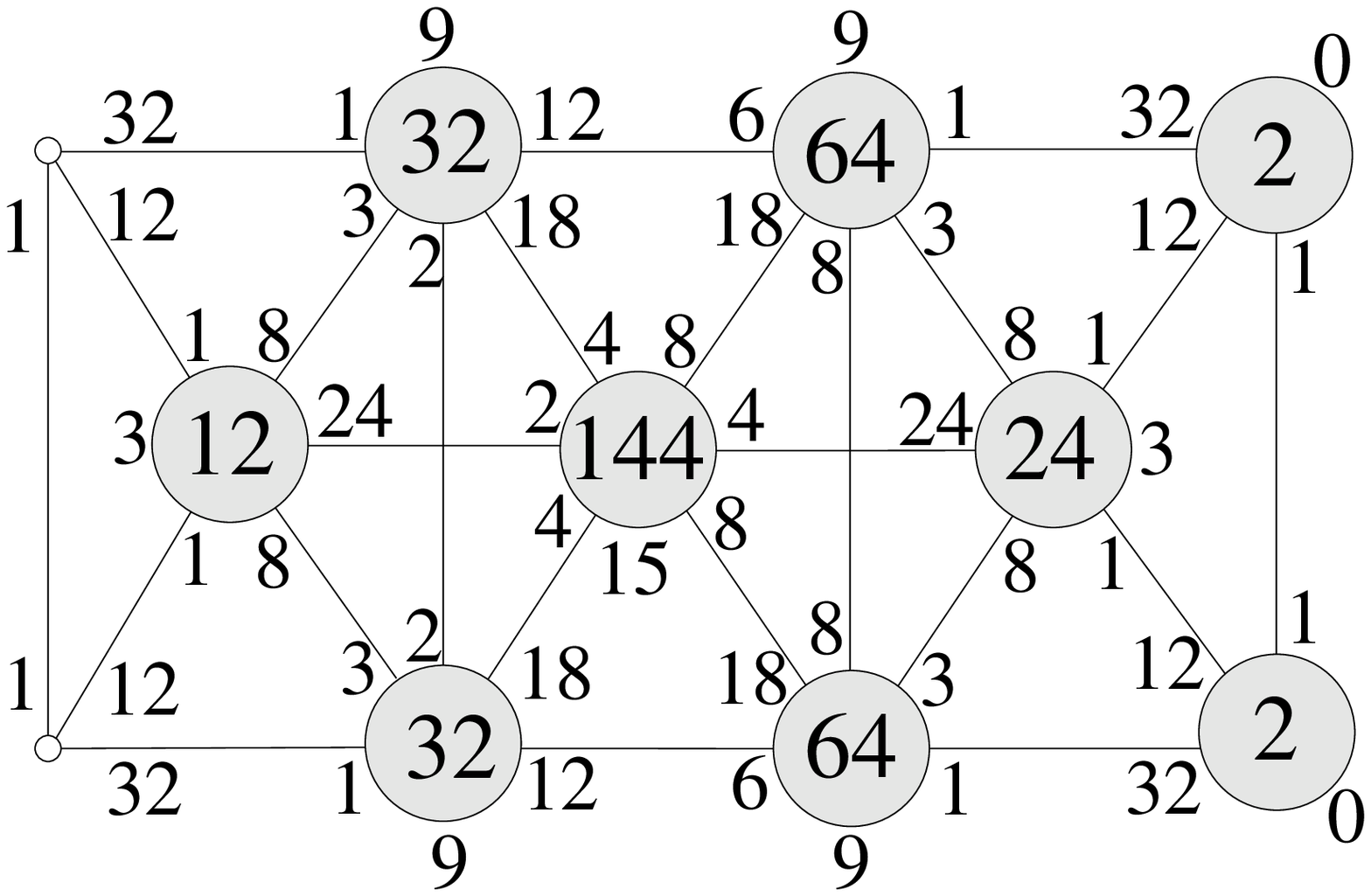,height=28mm}}

\vskip -1mm \noindent
\ \ \ \ \ {\footnotesize (a) AT4(1,2,3)}
 \ \ \ \ \ \ \ \ \ \ \ \ \ \ \ \ \ \ \ \ 
{\footnotesize (b) AT4(2,2,2)}
 \ \ \ \ \ \ \ \ \ \  \ \ \ \ \ \ \ \ \ 
{\footnotesize (c) AT4(4,2,2)}
 \ \ \ \ \ \ \ \ \ \ \ \ \ \ \ \ \ \ \
{\footnotesize (d) AT4(3,3,3)}

\bigskip
\baselineskip 9pt \par{\leftskip 1cm \rightskip 1cm \noindent
{\footnotesize 
Figure A.1: 1-homogeneous partition of (a) the Conway-Smith graph 
(b) the Johnson graph $J(8,4)$, \\
(c) the halved cube ${1\over 2}H(8,2)$, and (d) the $3.O_6^-(3)$.}
\par} \baselineskip=\normalbaselineskip

\medskip\medskip
\centerline{\psfig{figure=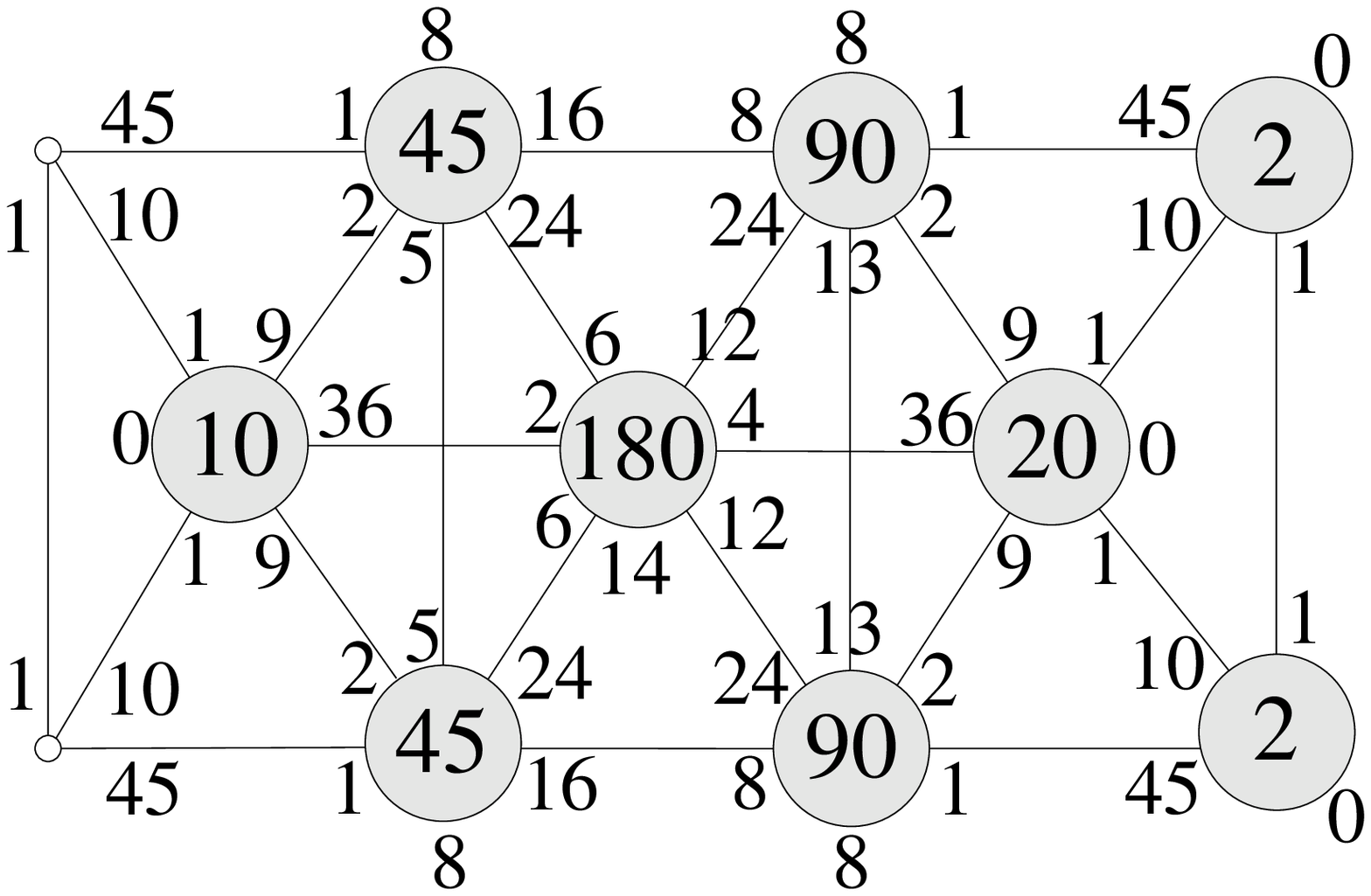,height=30mm}\ \ \ \ \ \ \ \ \ 
            \psfig{figure=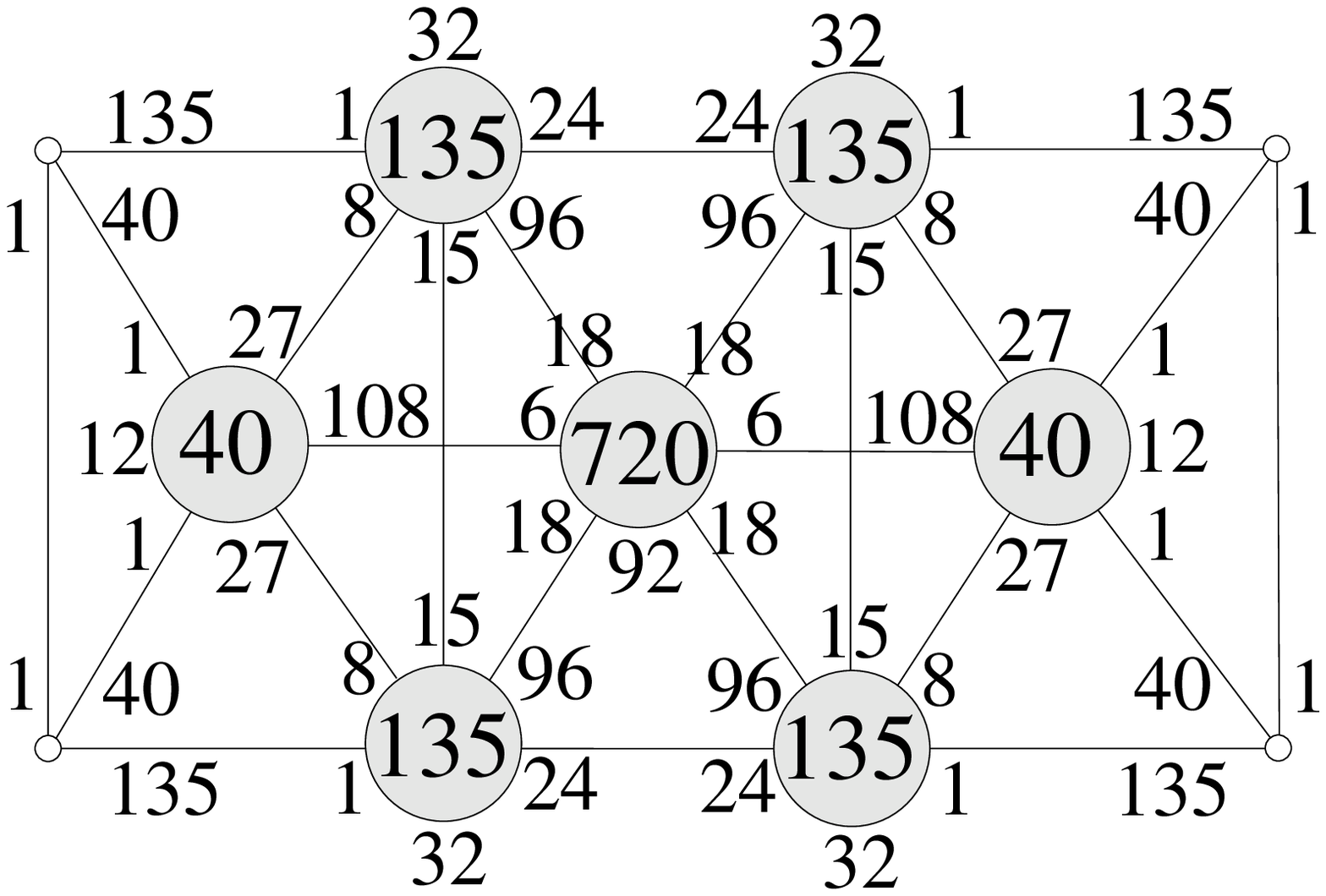,height=30mm}\ \ \ \ \ \ \ \ \
            \psfig{figure=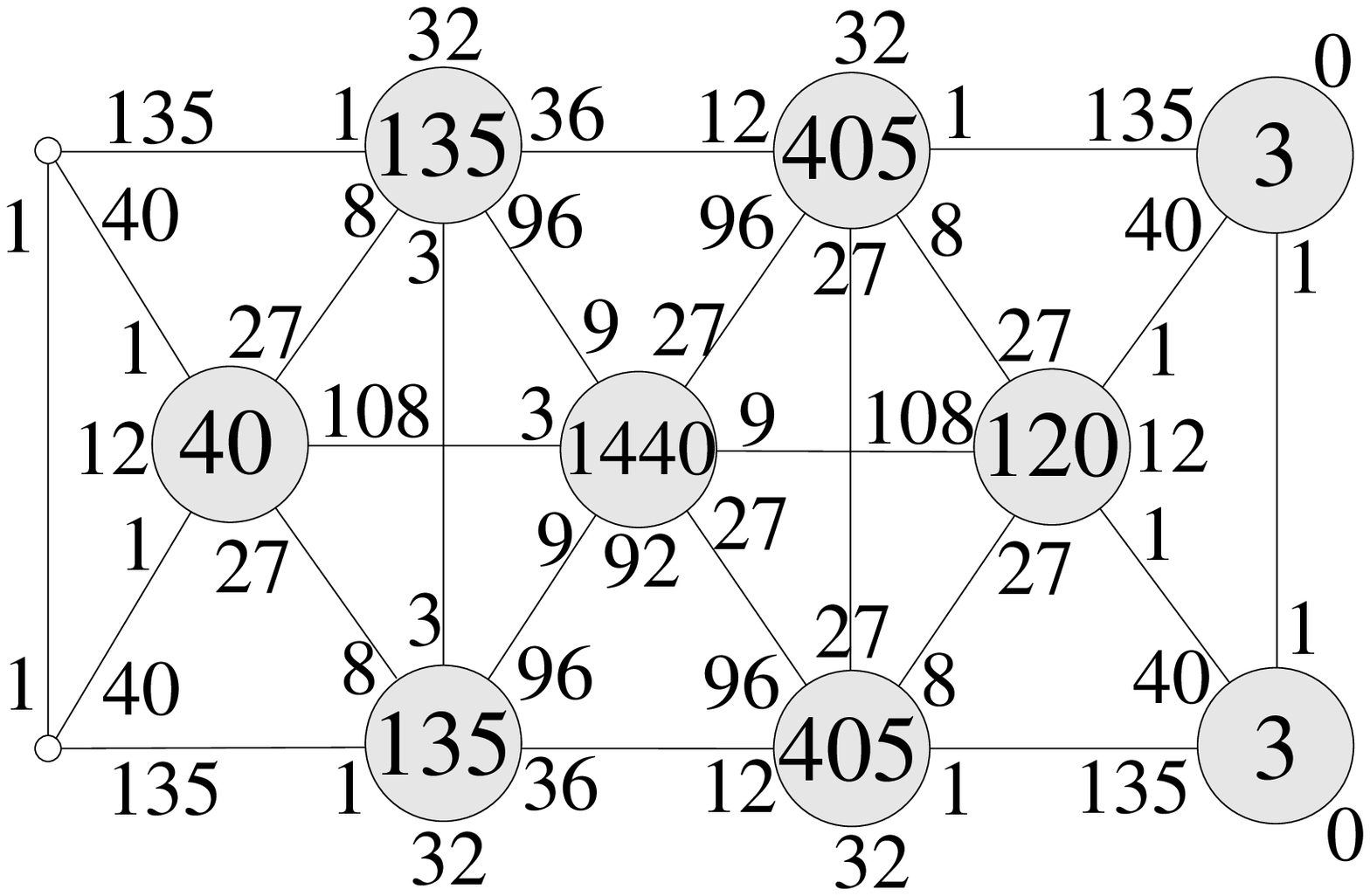,height=30mm}} 

\vskip -1mm \noindent
 \ \ \ \ \ \ \ \ \ \ \ \ \ \  
{\footnotesize (e) AT4(2,4,3)}
\ \ \ \ \ \ \ \ \ \ \ \ \ \ \ \ 
\ \ \ \ \ \ \ \ \ \ \
{\footnotesize (f) AT4(8,4,2)}
\ \ \ \ \ \ \ \ \ \ \ \ \ \ \ \ 
\ \ \ \ \ \ \ \ \ \ \ 
{\footnotesize (g) AT4(8,4,4)}

\medskip
\baselineskip 9pt \par{\leftskip 5mm \rightskip 5mm \noindent
{\footnotesize 
Figure A.2: 1-homogeneous partition of (e) the Soicher1 graph,  
(f) the Meixner1 graph, (g) the Meixner2 graph.}
\par} \baselineskip=\normalbaselineskip

\medskip\medskip
\centerline{\psfig{figure=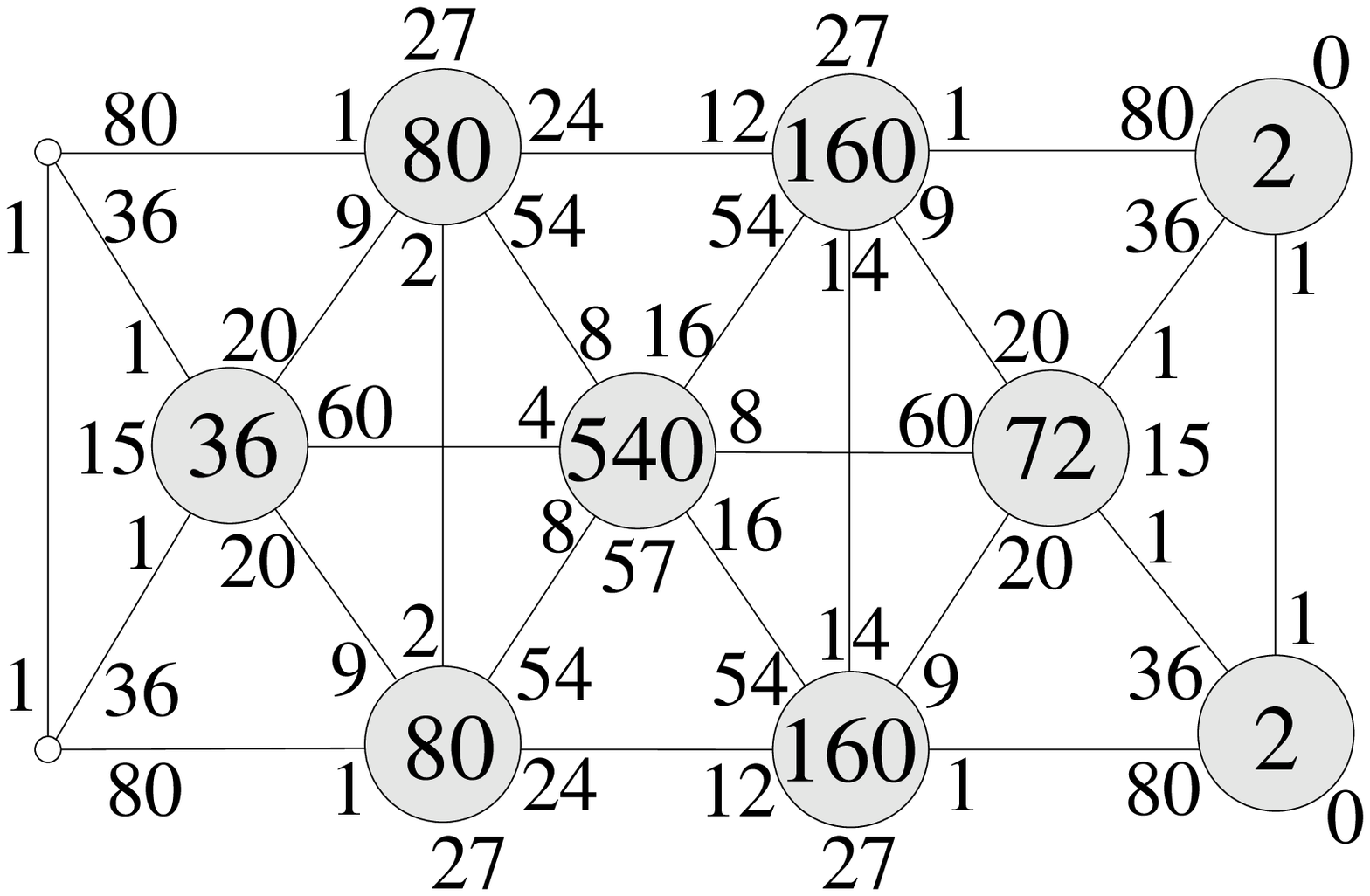,height=30mm} 
\ \ \ \ \ \ \ \ \ \psfig{figure=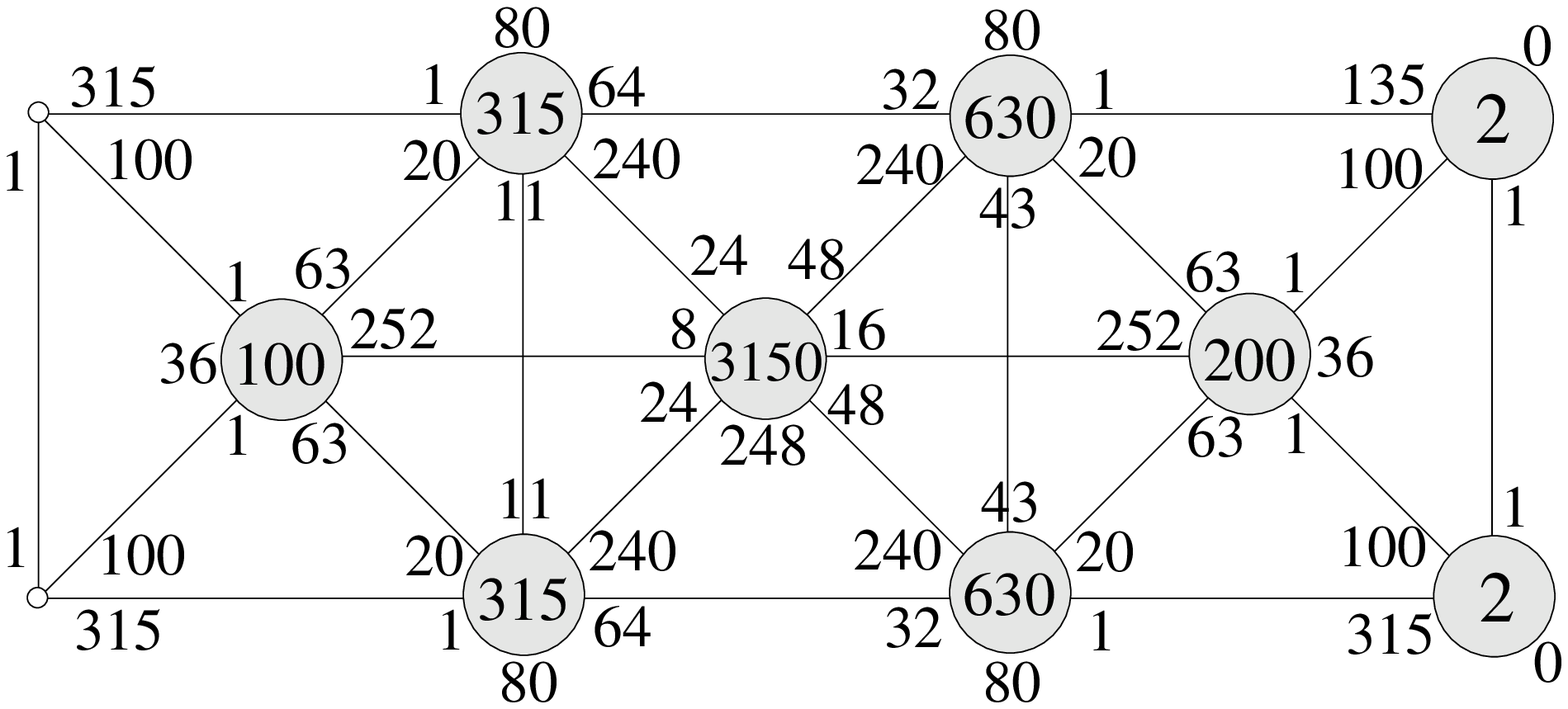,height=30mm}}

\vskip -1mm \noindent
\ \ \ \ \ \ \ \ \ \ \ \ \ \ \ \ \ \ \ \ \ \ \ \ \ \ \ 
{\footnotesize (h) AT4(9,3,3)}
\ \ \ \ \ \ \ \ \ \ \ \ \ \ \ \ \ \ \ \ \ \ \ \ \ \ \ \ \ \ \ \ 
{\footnotesize (i) AT4(20,4,3)}

{\footnotesize 
\medskip \centerline{\footnotesize
Figure A.3: 1-homogeneous partition of (h) the $3.O_7(3)$, 
(i) the Soicher2 graph.} 
\par} \baselineskip=\normalbaselineskip

\medskip\medskip
\centerline{\psfig{figure=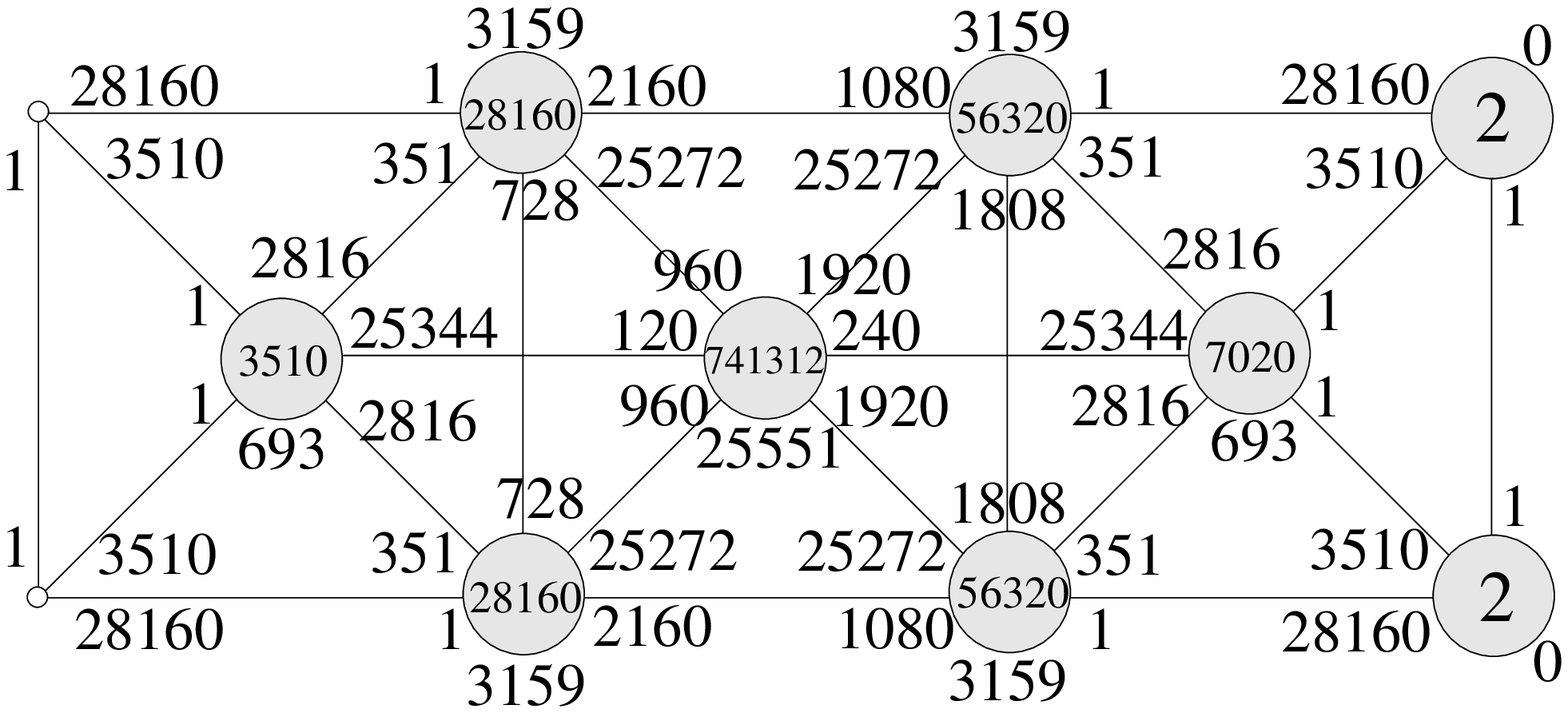,height=30mm}\ \ \ \ \ \ \ \ \ 
\ \ \ \psfig{figure=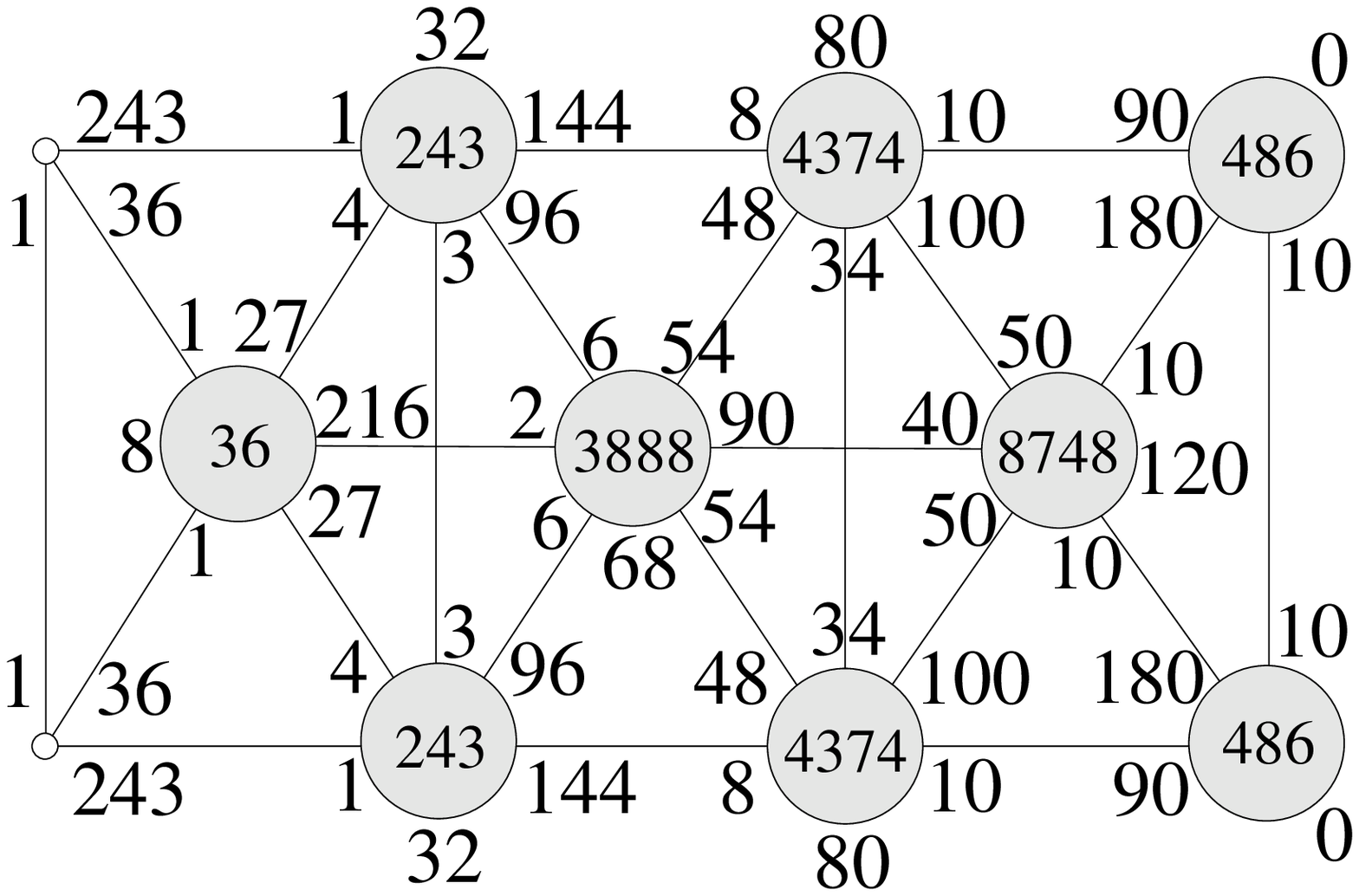,height=30mm}}

\vskip -1mm \noindent
\ \ \ \ \ \ \ \ \ \ \ \ \ \ \ \ \ \ \ \ \ \ \ \ \ 
\ \ \ \ \ \ \ \ \ {\footnotesize (j) AT4(351,9,3)}
\ \ \ \ \ \ \ \ \ \ \ \ \ \ \ \ \ \ \ \ \ \ \ \ \ \ \ \ \ \ \ \ \ \ \ \ 
\ \ \ \ \ \ \ \ \ \ {\footnotesize (k)}

\centerline{\footnotesize Figure A.4: 1-homogeneous partition of 
(j) the $3.Fi_{24}^-$ graph and (k) the Patterson graph.}

\end{document}